# Computational analysis on a linkage between generalized logit dynamic and discounted mean field game

Hidekazu Yoshioka[1,*]

[1]Japan Advanced Institute of Science and Technology, 1-1 Asahidai, Nomi, Ishikawa 923-1292, Japan

[*] Corresponding author: yoshih@jaist.ac.jp, ORCID:0000-0002-5293-3246

**Abstract**

Logit dynamics are dynamical systems describing transitions and equilibria of actions of interacting players under uncertainty. An uncertainty is embodied in logit dynamic as a softmax-type function often called a logit function originating from a maximization problem subjected to an entropic penalization. This study provides another explanation for the logit-type dynamic, particularly its logit function and player's heterogeneity, based on a discounted mean field game subjected to the costly decision-making of a representative player. A large-discount limit of the mean field game is argued to yield a logit dynamic. Further, mean field games that lead to classical and generalized logit dynamics are clarified and their well-posedness is discussed. Additionally, numerical methods based on a finite difference discretization for computing generalized logit dynamics and corresponding mean field games are presented. Numerical methods are applied to two problems arising in the management of resources and environment; one involves an inland fisheries management problem with legal and illegal anglers, while the other is a sustainable tourism problem. Particularly, cases that possibly lack the regularity condition to be satisfied for the unique existence of stationary solutions are computationally discussed.

## 1. Introduction

### 1.1 Background

Sustainable living with nature is a mission imposed on humans, for which the management of natural environment and associated ecosystem is a crucial issue. Particularly, natural resources including fish in a water body [1,2], groundwater in an aquifer [3,4], and timber in a forest [5,6] should be used without losing their renewability. Decision-making of and among individuals who benefit from using such resources plays a vital role in analyzing the sustainable development of human societies. Environmental- and resource-driven social dynamics, such as recreational fishing and ecotourism around geoparks, promote local sustainable development through outdoor education activities [7,8]. Game-theoretic models are pivotal tools for mathematically describing and analyzing social interactions among individuals and their consequences on the sustainable management of resources and the surrounding environment [9–11].

Evolutionary game theory is a branch of game theory applicable to a wide range of problems where social interaction among players follows a non-linear evolution equation (difference or differential equation) [12,13]. The core assumption in theory is that a player updates their current action in a utility-driven manner such that the player's action is determined through a dynamic comparison among their and other players' utilities. A solution to the governing evolution equation in an evolutionary game is often a probability measure or its density describing action profiles of players. Recent applications of evolutionary game models cover group formation [14], human behavior under injunctive social norms [15], cyber security [16], pollutant emission regulation [17], and low-carbon tourism [18].

An evolutionary game where players can choose their actions from a continuously parameterized set, arising in applications where they are allowed to use mixed strategies, associates a partial integro-differential equation (PIDE) as an evolution equation on a Banach space [19]. Examples include but are not limited to replicator dynamic [20,21], logit dynamic [22,23], and gradient dynamic [24,25], which use different protocols but share the common principle that resolution of an evolutionary game reduces to finding a proper solution to PIDE. Cheung [26] discussed some major evolutionary game models, such as replicator and logit dynamics, in a unified manner as a PIDE, called pairwise comparison dynamic (PCD). PCD implies that a broad class of PIDEs arising in evolutionary game theory is a special case of a unified non-local conservation law containing an inflow probability flux from other to current actions and an outflow probability flux from current to other action. PCD promotes PIDE studies in evolutionary games from the viewpoint of conservation laws, motivating researchers to seek for connections of PIDEs to dynamical systems arising in other theories. In this study, mean field game (MFG) is analyzed as an example among existing theories.

MFG is also a branch of game theory that is closer to optimal control than evolutionary game [27,28] and has recently been applied to a variety of social phenomena, including energy management for decarbonization [29], modeling prosumer behavior [30], and threshold-type resource harvesting [31]. MFG considers a population of interacting players as in evolutionary game but assumes a different mechanism of decision-making among them, such that players update their actions as feedback control variables considering the cumulated utility from current to future (often called value function). Due to a longer perspective than evolutionary games, the governing equation of an MFG model is given by both time-forward PIDE including the evolution equations as in evolutionary games and a time-backward PIDE to determine the optimal action based on a future prediction: an expectation conditioned on the current state. Numerical methods for MFGs for effective forward-backward solutions, such as the semi-Lagrangian and finite difference methods based on fictions play [32,33], higher-order discretization combined with a fixed-point iteration method [34], and kernel-based Monte Carlo method [35], have been proposed.

As aforementioned, actions of players in MFG are determined by forecasting based on current state, and its myopic limit is a PIDE in an evolutionary game. For a deterministic MFG model, Degond et al. [36] heuristically derived its myopic version based on a model predictive control. The method has been extended to MFG models, where diffusive stochastic differential equations (SDEs) drive players' actions [37] and pairwise logit dynamic as a jump-driven stochastic dynamical system [38]. Bertucci et al. [39] and Bardi and Cardaliaguet [40] studied MFGs subjected to discounted utilities and obtained rigorous convergence results under a large-discount limit, where a system of forward-backward PIDEs in an MFG gets closer to that of a PIDE in an evolutionary game. Thus, links between evolutionary game and MFG models have been studied in existing literature, but their implications to applied problems have limited exploration. Therefore, this study focuses on related evolutionary game and MFG models and their performance in environmental and resource management, where game-theoretic models play a leading role.

### 1.2 Contribution and objective

The objective of this study includes:

- ✓ To investigate an evolutionary game model that covers conventional logit dynamics and derive its MFG counterpart, and
- ✓ To compare the evolutionary game and MFG models through their applications to management problems on resource and environment.

The contributions of this study to achieve the objectives cover both theory and application of evolutionary game and MFG models.

To achieve the first objective, a classical logit dynamic [22] of utility-maximizing players is generalized in a continuous action space as a special case of PCDs through an extension of the logit function (softmax part in the dynamic). The generalization is different from those examined in existing literature [23,41] in that the generalization of this study has a different form of non-locality and the difference is due to the MFG counterpart. Moreover, the logit dynamic of this study assumes the existence of multiple groups [e.g., 42,43] to account for discrete heterogenous preferences among players. The well-posedness of the evolution equation, called generalized logit dynamic (GLD), is studied, and its solution is found to be a time-dependent probability measure in a Banach space equipped with a total variation norm.

This paper presents an MFG model with a discount, whose objective function to be maximized by a representative player contains the utility of the GLD scaled by a constant parameter and cost to update their actions based on the latest available information. The discount rate plays a vital role in MFG. This study heuristically demonstrates that designing a proper cost term leads to a forward-backward PIDE system that formally has the forward part as a GLD, where the utility is replaced by a value function. Further, taking a large discount rate is demonstrated to reduce the time-backward PIDE in MFG to an algebraic equation and the time-forward PIDE to GLD. The optimal action in MFG arises as the logit function in the GLD, providing its novel mean field description. Particularly, the first objective is new even for the classical logit case [41] because the MFG counterpart has not been discussed in literature to the best of the author's knowledge.

To achieve the second objective, finite difference methods for computing GLD and MFG are discussed, and the forward-backward nature of the latter is computed by a fixed-point iteration with a relaxation. Two application problems related to resource and environmental management are discussed, including sustainable tourism and recreational fishing by legal and illegal anglers. The applications focused on in this study are categorized as potential games [44], where the utility is possibly concave (often called monotone). This property often guarantees unique existence and stability of solutions to the PIDEs in game theory [45,46]. A potential game is a well-studied subject in literature, but its understanding from the perspective of inland fisheries management is limited. This study focuses on novel, simple, and less mechanistic behavioral modeling of anglers and tourists than the microscopic models discussed in existing literature [47,48]. The potential game description allows effective behavioral modeling of the anglers, showing better suitability to consistent modeling between evolutionary game and MFG models.

The remainder of this paper is organized as follows. **Section 2** presents preliminaries, introduces the GLD, and reconsiders it as a PCD. **Section 3** analyzes an MFG whose formal limit is a GLD and presents few specific examples, where the source function of the MFG is found explicitly. **Section 4** is devoted to computational analysis of the GLD and MFG, focusing on applied problems. **Section 5** concludes this study and presents the perspectives. Auxiliary results including proofs of propositions are placed in **Appendix**.

## 2. Generalized logit dynamic

### 2.1 Preliminaries

For this study, preliminaries present several notions used based on literature [e.g., 19]. Time $t$ is a non-negative continuous parameter. The space of continuous actions of players is given by the unit compact interval $\Omega = [0,1]$. The collection of all $\sigma$-algebras on $\Omega$ is denoted as $\mathcal{B}$, and the space of all finite signed measures in $(\mathcal{B},\Omega)$ is given as $\mathcal{M}(\mathcal{B},\Omega)$. The space $\mathcal{M}(\mathcal{B},\Omega)$ is equipped with the total variation norm given for any $\mu \in \mathcal{M}(\mathcal{B},\Omega)$ as $\|\mu\| = \sup_g \left| \int_0^{+\infty} g(x)\mu(\mathrm{d}x) \right|$. Herein, the supremum is taken with respect to all measurable functions satisfying $\sup_{x\in\Omega} |g(x)| \le 1$. A set of non-negative measures with the total mass $m > 0$ is defined as

$$\mathcal{P}_m(\mathcal{B},\Omega) = \left\{ \mu \in \mathcal{M} \,\middle|\, \mu(A) \ge 0 \text{ for any } A \in \mathcal{B},\ \mu(\Omega) = m \right\}. \tag{1}$$

A set of signed measures having the norm not larger than 2 is defined as $\mathcal{M}_2(\mathcal{B},\Omega)=\{\mu\in\mathcal{M} \mid \|\mu\|\leq 2\}$. The space $\mathcal{P}_m$ serves as a basis of players' actions, while $\mathcal{M}_2$ is introduced for technical reasons to prove unique existence of solutions to the GLD (**Proof of Proposition 2**).

Additionally, notions of types of players are introduced. Players of distinct types have different utility functions with each other. The type of players is represented by a natural number in the discrete set $\Xi_I=\{1,2,3,...,I\}$ with some $I\in\mathbb{N}$. Thus, the type of players is assumed to be fixed, while their actions are updated in time. Then, a set of probability measures that concurrently represent players' actions and types is introduced. For each $I\in\mathbb{N}$, the probability masses $\{m_i\}_{i=1,2,3,...,I}$ are introduced: $m_i>0$ and $\sum_{i=1}^{I} m_i=1$. The series $\{m_i\}_{i=1,2,3,...,I}$ represents proportion of players of each type. The distribution of players' actions of the $i$ th type is represented by a measure $\mu_i\in\mathcal{P}_{m_i}(\mathcal{B},\Omega)\subset\mathcal{M}_2$. Thus, $\sum_{i=1}^{I}\mu_i(\Omega)=\sum_{i=1}^{I} m_i=1$.

In the sequel, the arguments $(\mathcal{B},\Omega)$ are omitted to simplify notations from $\mathcal{M}(\mathcal{B},\Omega)$ as $\mathcal{M}$. The same applies to $\mathcal{P}_m(\mathcal{B},\Omega)$ land $\mathcal{M}_2(\mathcal{B},\Omega)$. Then, $I\in\mathbb{N}$ is fixed for the rest of this study unless otherwise specified. Therefore, the subscript $I$ is omitted from $\Xi_I$. If a measure $\mu$ is time-dependent, it can be written as $\mu(t)$ at time $t$ and its integration on $A\in\mathcal{B}$ as $\mu(t,A)$.

## 2.2 Utility and deformed exponential

The utility $U$ in this study satisfies **Assumption 1** below unless otherwise specified. Hereafter, $\{\mu\}$ is an abbreviation of the collection $\{\mu_i\}_{i=1,2,3,...,I}$.

***Assumption 1*** *Utility* $U:\Xi\times\Omega\times(\mathcal{M})^I\to\mathbb{R}$ *satisfies the following conditions; there exist constants* $L>0$ *and* $d\in(0,1]$ *such that*

$$\left|U_i(x,\{\mu\})\right|\leq L \quad \text{for all} \quad (i,x,\{\mu\})\in\Xi\times\Omega\times(\mathcal{M}_2)^I, \tag{2}$$

$$\left|U_i(x,\{\mu\})-U_i(y,\{\mu\})\right|\leq L|x-y|^d \quad \text{for all} \quad (i,x,y,\{\mu\})\in\Xi\times\Omega\times\Omega\times(\mathcal{M}_2)^I, \tag{3}$$

$$\left|U_i(x,\{\mu\})-U_i(x,\{\nu\})\right|\leq L\sum_{j=1}^{I}\left\|\mu_j-\nu_j\right\| \quad \text{for all} \quad (i,x,\{\mu\},\{\nu\})\in\Xi\times\Omega\times(\mathcal{M}_2)^I\times(\mathcal{M}_2)^I. \tag{4}$$

Herein, condition (2) bounds the utility, while conditions (3) and (4) guarantee its continuity. The Hölder and Lipschitz continuity conditions are required to establish the unique existence of a solution to the GLD. These conditions are common in the logit dynamics [22] and is therefore not restrictive.

A deformed exponential function is used as a generalized exponential function, and the Tsallis one is considered due to its wide applications [49,50]:

$$\exp_q(z)=\begin{cases}(1+(1-q)z)_+^{\frac{1}{1-q}} & (q>1 \text{ and } 1+(1-q)z>0, \text{ or } 0<q<1 \text{ and } z\in\mathbb{R}) \\ +\infty & (q>1 \text{ and } 1+(1-q)z\le 0) \\ \exp(z) & (q=1,\ z\in\mathbb{R})\end{cases}, \tag{5}$$

where $(z)_+=\max\{0,z\}$ ($z\in\mathbb{R}$). The deformed exponential is not well-defined for the second case of (5), although it is continuous and increasing for the other cases. The cases $q>1$, $q=1$, and $0<q<1$ correspond to less sensitive, conventional, and more sensitive softmax functions.

### 2.3 Formulation and well-posedness

For each $i\in\Xi$ and $t>0$, the GLD is formulated as

$$\underbrace{\frac{\mathrm{d}}{\mathrm{d}t}\mu_i(t,A)}_{\text{Temporal change}}=\underbrace{m_i\mathfrak{L}_i\mu(t,A)}_{\text{Inflow from outside } A}-\underbrace{\theta_i(q,\eta)\mu_i(t,A)}_{\text{Outflow from } A} \quad \text{for any } A\in\mathcal{B} \tag{6}$$

with

$$\mathfrak{L}_i\mu(t,A)=\int_{y\in A}\int_{w\in\Omega}\left\{\frac{\exp_q\left(\eta^{-1}\Delta U_i(y,w,\{\mu\})\right)}{\int_\Omega \exp_q\left(\eta^{-1}\Delta U_i(z,w,\{\mu\})\right)\mathrm{d}z}\right\}^q \mathrm{d}w\mathrm{d}y \tag{7}$$

$$\Delta U_i(x,y,\{\mu\})=U_i(x,\{\mu\})-U_i(y,\{\mu\}), \tag{8}$$

and

$$\theta_i(q,\eta)=\int_{y\in\Omega}\int_{w\in\Omega}\left\{\frac{\exp_q\left(\eta^{-1}\Delta U_i(w,y,\{\mu\})\right)}{\int_\Omega \exp_q\left(\eta^{-1}\Delta U_i(z,y,\{\mu\})\right)\mathrm{d}z}\right\}^q \mathrm{d}w\mathrm{d}y>0. \tag{9}$$

The PIDE (6) is subject to an initial condition $\mu_i(0)\in\mathcal{P}_{m_i}$ ($i\in\Xi$; recall that $\sum_{i=1}^{I}m_i=1$). Herein, $\eta>0$ is an uncertainty parameter such that a larger value represents larger uncertainty that each player faces, and hence their distribution of actions can become more blurred. The time derivative in (6) is defined in a strong sense (e.g., Definition 1.5 in [51]):

$$\lim_{\varepsilon\to 0}\left\|\frac{\mathrm{d}}{\mathrm{d}t}\mu_i(t,A)-\frac{\mu_i(t+\varepsilon,A)-\mu_i(t,A)}{\varepsilon}\right\|=0. \tag{10}$$

Before analyzing the well-posedness of GLD (6), the difference between the proposed and existing logit dynamics requires investigation. First, as in the existing logit dynamics, the right-hand side of (6) is formally represented by (an integration of) a measure minus another measure. The GLD reduces to the classical logit dynamic [22] when $I=1$ and $q=1$. The case $I>1$ and $q=1$ is like the dynamic discussed by Lahkar and Ramani [45]. However, if $q\neq 1$, such a relationship does not hold true for logit dynamics based on deformed exponentials discussed by Lahkar et al. [41] and Yoshioka [23]. Indeed, the logit function $\mathfrak{L}_i\mu$ in the context of exiting literature is given by

$$\mathfrak{L}_i\mu(t,A)=\int_A \frac{\exp_q\left(\eta^{-1}U\left(y,\{\mu\}\right)\right)}{\int_\Omega \exp_q\left(\eta^{-1}U_i\left(z,\{\mu\}\right)\right)\mathrm{d}z}\mathrm{d}y\,. \tag{11}$$

Our dynamic has a utility difference inside a deformed exponential, while those by Lahkar et al. [41] and Yoshioka [23] have a utility. The difference is due to the consideration of uncertainty in measuring a utility difference by logit dynamic, as inferred in **Section 3**, while those in existing literature is due to that in measuring the utility itself. The qualitative difference vanishes for the classical case ( $q=1$ ) since

$$\int_{y\in A}\int_{w\in\Omega}\left\{\frac{\exp_1\left(\eta^{-1}\Delta U_i\left(y,w,\{\mu\}\right)\right)}{\int_\Omega \exp_1\left(\eta^{-1}\Delta U_i\left(z,w,\{\mu\}\right)\right)\mathrm{d}z}\right\}^1 \mathrm{d}w\mathrm{d}y=\int_{y\in A}\frac{\exp_1\left(\eta^{-1}\Delta U_i\left(y,\{\mu\}\right)\right)}{\int_\Omega \exp_1\left(\eta^{-1}\Delta U_i\left(z,\{\mu\}\right)\right)\mathrm{d}z}\mathrm{d}y\,. \tag{12}$$

To well-pose the GLD (6), the following condition is assumed for the rest of this paper.

***Assumption 2*** *If* $q\neq 1$, *then there exists a constant* $\vartheta>0$ *such that*

$$1-2\left|1-q\right|\eta^{-1}L\geq\vartheta\,. \tag{13}$$

This ensures the well-refinedness of $\mathfrak{L}_i\mu$ ; if $q\neq 1$ and (13) is satisfied, then

$$1+\left(1-q\right)\eta^{-1}\Delta U_i\left(x,y,\{\mu\}\right)\geq 1-\left|1-q\right|\eta^{-1}\times 2L\geq\vartheta>0\,, \tag{14}$$

and the second case in (5) does not occur, and moreover we have that both the denominator and numerator of $\mathfrak{L}_i\mu$ are positive, and it is therefore well-defined. **Assumption 2** implies that considering a small-uncertainty limit $\eta\to+0$ may be difficult if $q\neq 1$, which is therefore a theoretical restriction, but our computational results demonstrate that **Assumption 2** can be removed in some case if $0<q<1$.

***Proposition 1***

*The GLD (6) admits a unique solution* $\{\mu(t)\}=\{\mu_i(t)\}_{i=1,2,3,\ldots,I}$ *such that* $\mu_i(t)\in\mathcal{P}_{m_i}$ *for all* $i\in\Xi$.

**Proposition 1** is proven in the following way. First, the Lipschitz continuity of the right-hand side (7) is shown for $\{\mu(t)\}\in(\mathcal{M}_2)^I$. Then, a modified GLD is considered; for each $i\in\Xi$ and $t>0$,

$$\frac{\mathrm{d}}{\mathrm{d}t}\mu_i(t,A)=\left(2-m_i^{-1}\left\|\mu_i(t,\cdot)\right\|\right)_+\left(m_i\mathfrak{L}_i\mu(t,A)-\theta_i(q,\eta)\mu_i(t,A)\right)\ \text{ for any }\ A\in\mathcal{B}\,. \tag{15}$$

Utilizing the coefficient $\left(2-m_i^{-1}\left\|\mu_i(t,\cdot)\right\|\right)_+$, the right-hand side of (15) is proven to be Lipschitz continuous with respect to $\{\mu(t)\}\in(\mathcal{M})^I$. Then, the generalized Picard–Lindelöf theorem (p. 79 [52]) is applied to the modified GLD, confirming a unique solution $\{\mu(t)\}\in(\mathcal{M})^I$ ( $t\geq 0$ ). The generalized

Picard–Lindelöf theorem is not applicable to the case where $(\mathcal{M})^I$ is replaced by $(\mathcal{M}_2)^I$ as it *a priori* requires $\{\mu(t)\} \in (\mathcal{M}_2)^I$, which is not guaranteed at this stage. Thereafter, the solution is demonstrated to satisfy $\mu_i(t) \in \mathcal{P}_{m_i}$ ($i \in \Xi$) to conclude $\left(2 - m_i^{-1}\|\mu_i(t,\cdot)\|\right)_+ = 1$. Then, the Lipschitz continuity of the right-hand side of (7) proven in the first step guarantees that the solution is the desired one. Cheung [26] also used a qualitatively the same strategy of proof to show the well-posedness of a pair-wise comparison dynamic; however, the problem of the form proposed in this paper has not been discussed, which should be checked by analytical calculations.

### 2.4 Generalized logit dynamic as a pairwise comparison dynamic

A PCD is a general class of evolutionary game model [e.g., 26]. In our context, for each $t > 0$ and $i \in \Xi$, PCD reads

$$\begin{aligned}\underbrace{\frac{\mathrm{d}}{\mathrm{d}t}\mu_i(t,A)}_{\text{Temporal change on }A} &= \underbrace{\int_{x\in A}\int_{y\in\Omega}\rho_i(x,y,\{\mu\})\mu_i(t,\mathrm{d}y)\mathrm{d}x}_{\text{Inflow from outside }A} \\ &\quad -\underbrace{\int_{x\in A}\int_{y\in\Omega}\rho_i(y,x,\{\mu\})\mathrm{d}y\mu_i(t,\mathrm{d}x)}_{\text{Outflow from }A}\end{aligned} \quad \text{for any } A \in \mathcal{B}, \tag{16}$$

where $\rho_i(x,y,\{\mu\})$ with $\rho: \Omega\times\Omega\times(\mathcal{M})^I \to [0,+\infty)$ represents the transition rate of players' actions for the $i$ th type from $y$ to $x$. PCD is a generalization of major evolutionary game models including logit and replicator dynamics. The GLD of this study is not an exception. Indeed, set

$$\rho_i(x,y,\{\mu\}) = \int_{w\in\Omega}\left\{\frac{\exp_q\left(\eta^{-1}\Delta U_i(x,w,\{\mu\})\right)}{\int_\Omega \exp_q\left(\eta^{-1}\Delta U_i(z,w,\{\mu\})\right)\mathrm{d}z}\right\}^q \mathrm{d}w \quad \text{(formally independent from } y\text{)}. \tag{17}$$

Then, $\int_{y\in\Omega}\mu_i(\mathrm{d}y) = m_i$ is recalled to obtain

$$\begin{aligned}\int_{x\in A}\int_{y\in\Omega}\rho_i(x,y,\{\mu\})\mu_i(t,\mathrm{d}y)\mathrm{d}x &= \int_{x\in A}\int_{y\in\Omega}\int_{w\in\Omega}\left\{\frac{\exp_q\left(\eta^{-1}\Delta U_i(x,w,\{\mu\})\right)}{\int_\Omega \exp_q\left(\eta^{-1}\Delta U_i(z,w,\{\mu\})\right)\mathrm{d}z}\right\}^q \mathrm{d}w\mu_i(\mathrm{d}y)\mathrm{d}x \\ &= \int_{x\in A}\int_{w\in\Omega}\left\{\frac{\exp_q\left(\eta^{-1}\Delta U_i(x,w,\{\mu\})\right)}{\int_\Omega \exp_q\left(\eta^{-1}\Delta U_i(z,w,\{\mu\})\right)\mathrm{d}z}\right\}^q \mathrm{d}w\mathrm{d}x\left(\int_{y\in\Omega}\mu_i(\mathrm{d}y)\right) \\ &= m_i\mathfrak{L}_i\mu(t,A)\end{aligned} \tag{18}$$

and

$$
\begin{aligned}
&\int_{x\in A}\int_{y\in\Omega}\rho_i\left(y,x,\{\mu\}\right)\mathrm{d}y\mu_i\left(t,\mathrm{d}x\right)\\
&=\int_{x\in A}\int_{y\in\Omega}\int_{w\in\Omega}\left\{\frac{\exp_q\left(\eta^{-1}\Delta U_i\left(w,y,\{\mu\}\right)\right)}{\int_\Omega \exp_q\left(\eta^{-1}\Delta U_i\left(z,y,\{\mu\}\right)\right)\mathrm{d}z}\right\}^q \mathrm{d}w\mathrm{d}y\mu_i\left(\mathrm{d}x\right)\\
&=\int_{x\in A}\mu_i\left(\mathrm{d}x\right)\int_{y\in\Omega}\int_{w\in\Omega}\left\{\frac{\exp_q\left(\eta^{-1}\Delta U_i\left(w,y,\{\mu\}\right)\right)}{\int_\Omega \exp_q\left(\eta^{-1}\Delta U_i\left(z,y,\{\mu\}\right)\right)\mathrm{d}z}\right\}^q \mathrm{d}w\mathrm{d}y\\
&=\theta_i\left(q,\eta\right)\mu_i\left(A\right)
\end{aligned}
\tag{19}
$$

Herein, (17)–(19) prove that the GLD of this study is a special case of PCD. Due to PCD representations, the mass conservation of the GLD follows immediately:

$$
\begin{aligned}
\frac{\mathrm{d}}{\mathrm{d}t}\mu_i\left(t,\Omega\right)&=\int_{x\in\Omega}\int_{y\in\Omega}\rho_i\left(x,y,\{\mu\}\right)\mu_i\left(t,\mathrm{d}y\right)\mathrm{d}x-\int_{x\in\Omega}\int_{y\in\Omega}\rho_i\left(y,x,\{\mu\}\right)\mathrm{d}y\mu_i\left(t,\mathrm{d}x\right)\\
&=\int_{x\in\Omega}\int_{y\in\Omega}\rho_i\left(x,y,\{\mu\}\right)\mu_i\left(t,\mathrm{d}y\right)\mathrm{d}x-\int_{x\in\Omega}\int_{y\in\Omega}\rho_i\left(x,y,\{\mu\}\right)\mathrm{d}x\mu_i\left(t,\mathrm{d}y\right).\\
&=0
\end{aligned}
\tag{20}
$$

***Remark 1*** If $0<q<1$, then formally

$$
\rho_i\left(x,y,\{\mu\}\right)=\int_{w\in\Omega}\left\{\frac{\exp_q\left(\eta^{-1}\Delta U_i\left(x,w,\{\mu\}\right)\right)}{\int_\Omega \exp_q\left(\eta^{-1}\Delta U_i\left(z,w,\{\mu\}\right)\right)\mathrm{d}z}\right\}^q \mathrm{d}w\to\int_{w\in\Omega}\left\{\frac{\left(\Delta U_i\left(x,w,\{\mu\}\right)\right)_+^{\frac{1}{1-q}}}{\int_\Omega\left(\Delta U_i\left(z,w,\{\mu\}\right)\right)_+^{\frac{1}{1-q}}\mathrm{d}z}\right\}^q \mathrm{d}w \tag{21}
$$

under $\eta\to+0$. Thus, GLD formally reduces to PCD (16) with $\rho_i$ in the right-most side of (21), implying that the vanishing-uncertainty limit is in general different from the Nash equilibria as maximizers of the utility $U_i$, because solutions to GLD should depend on $q$ that is external to the utility [23]. Therefore, a particular property is equipped with GLD such that $0<q<1$, representing heavier and more persistent uncertainties in decision-making of players than that in the case of $q=1$.

## 3. Connection to mean field description

### 3.1 Mean field game

This study heuristically formulates an MFG model that is connected to the GLD. Assuming that there exists $N\in\mathbb{N}$ players, the $n$ th player has some fixed type $i\in\Xi$, whose action is represented by a process $X^{(n,i)}=\left(X_t^{(n,i)}\right)_{t\geq 0}$. The empirical probability measure $\mu^{(N,i)}$ of actions of type-$i$ players formally reads

$$
\mu^{(N,i)}\left(t,\mathrm{d}x\right)=\frac{1}{N^{(i)}}\sum_{n=1}^{N^{(i)}}\delta\left(x-X_t^{(n,i)}\right),\quad t\geq 0, \tag{22}
$$

where $\delta(\cdot)$ is the Dirac delta and $N^{(i)}\in\mathbb{N}$ is the total number of type-$i$ players. Thus, $\sum_{i=1}^{I}N^{(i)}=N$.

The process $X^{(n,i)}$ is assumed to follow the jump SDE given by

$$
\mathrm{d}X_t^{(n,i)}=\int_0^1\left(z-X_{t-}^{(n,i)}\right)P^{(n,i)}\left(\mathrm{d}t,\mathrm{d}z\right),\quad t>0, \tag{23}
$$

subjected to an initial condition $X_0^{(n,i)} \in \Omega$, where $X_{t-}^{(n,i)}$ is the left limit of $X_t^{(n,i)}$. Herein, $P^{(n,i)}$ is a Poisson random measure on $(0,+\infty)\times\Omega$, whose compensated measure is $\left\{\varphi_t^{(n,i)}(z)\right\}^q \mathrm{d}z\mathrm{d}t$ with $\varphi_t^{(n,i)}$ being a distortion of the uniform probability measure $\mathrm{d}z$ on $\Omega$ with $\int_\Omega \varphi_t^{(n,i)}(z)\mathrm{d}z = 1$, which is a random field that serves as a control variable in the MFG. The appearance of $\left\{\varphi_t^{(n,i)}(z)\right\}^q$ in the compensated measure is due to the consideration of a uniform distribution distortion in an escort sense, which is consistent with the generalized divergence [53,54] and leads to a closed-form optimal control, as discussed further. Moreover, the formalism reduces to the classical distortion formulation of the compensator if $q=1$ [55]. Hence, the proposed formalism generalizes it. Elements in $\left\{P^{(n)}\right\}_{n=1,2,3,\ldots,N}$ are assumed to be mutually independent.

For each $n \in \{1,2,3,\ldots,N\}$ and sequence $\{x_k\}_{k=1,2,3,\ldots,N}$, set $x_{-n} = \{x_k\}_{\substack{k=1,2,3,\ldots,N \\ k\neq n}}$. Each player updates his/her action such that the objective function (26) is maximized by optimizing $\left(\varphi_t^{(n,i)}(\cdot)\right)_{t>0}$:

$$J^{(n,i)}(t,x_n,x_{-n}) = \mathbb{E}\left[\left.\begin{array}{l} \underbrace{w\int_t^T\int_\Omega e^{-\delta(t-s)}U_i\left(X_s^{(n,i)},\left\{\mu^{(N,i)}(s,\cdot)\right\}_{i=1,2,3,\ldots,I}\right)\mathrm{d}z\mathrm{d}s}_{\text{Cumulated utility}} \\ -\underbrace{\int_t^T\int_\Omega e^{-\delta(t-s)}\phi\left(\varphi_s^{(n,i)}(z)\right)\mathrm{d}z\mathrm{d}s}_{\text{Cumulated control cost}} \end{array}\right|\left\{X_t^{(n,i)} = x_n\right\}_{n=1,2,3,\ldots,N}\right] \tag{24}$$

for each $t\in(0,T)$ and $(x_n,x_{-n})\in\Omega^N$, where $T>0$ is a fixed terminal time, $\delta>0$ is a discount rate, $w>0$ is a weight constant, and $\phi:[0,+\infty)\to[0,+\infty)$ is a function representing control cost designed in the later sections. The discount rate $\delta$ modulates perspectives of players such that they are more myopic (less myopic) for a larger value (smaller value) of $\delta$. The objective function (24) is an expectation conditioned on the current state $(t,x_n,x_{-n})$ of players and contains the utility to be maximized (first term in the right-hand side) and the cost of controlling their action distributions. Intuitively, the cost should be minimal for a uniform distribution on $\Omega$ that should be the most randomized one is required. The terminal time $T$ is introduced for a technical reason, inferring that the objective function (24) can be reasonably considered as approximation of an infinite-horizon one at time $t<<T$ provided that $T$ is sufficiently large: the turnpike. The discount rate $\delta$ connects the MFG with the GLD, as discussed further.

A large population limit is considered, where the convergence $\left\{N^{-1}N^{(i)}\right\}_{i=1,2,3,\ldots,I} \to \{m_i\}_{i=1,2,3,\ldots,I}$ holds under $N\to+\infty$. By assuming that players' actions are homogenous for each $i$ th type, the SDE governing the type-$i$ representative player's action $X^{(i)} = \left(X_t^{(i)}\right)_{t\geq 0}$ is inferred, with abusing notations,

$$\mathrm{d}X_t^{(i)} = \int_0^1\left(z - X_{t-}^{(i)}\right)P^{(i)}(\mathrm{d}t,\mathrm{d}z),\quad t>0, \tag{25}$$

where $P^{(i)}$ is a Poisson random measure on $(0,+\infty)\times\Omega$, whose compensated measure is $\left\{\varphi_t^{(i)}(z)\right\}^q \mathrm{d}z\mathrm{d}t$ with a random field $\varphi_t^{(i)}$ being a distortion of the uniform probability measure $\mathrm{d}z$ on $\Omega$, satisfying $\int_\Omega \varphi_t^{(i)}(z)\mathrm{d}z = 1$. Then, the mean-field counterpart of the objective function (24) can be inferred as

$$J^{(i)}(t,x) = \mathbb{E}\left[ w\int_t^T \int_\Omega e^{-\delta(t-s)} U_i\left(X_s^{(i)},\{\mu(s,\cdot)\}\right)\mathrm{d}z\mathrm{d}s - \int_t^T \int_\Omega e^{-\delta(t-s)} \phi\left(\varphi_s^{(i)}(z)\right)\mathrm{d}z\mathrm{d}s \middle| X_t^{(i)} = x \right] \quad (26)$$

and its maximized version, the value function, as

$$\Phi^{(i)}(t,x) = \sup_{\varphi^{(i)}} J^{(i)}(t,x), \quad t\in[0,T] \text{ and } x\in\Omega. \quad (27)$$

Furthermore, its time-backward PIDE (the Hamilton–Jacobi–Bellman (HJB) equation associated with the MFG) for $i\in\Xi$ is given by (MFG with jumps, [56–59])

$$-\frac{\partial\Phi^{(i)}(t,x)}{\partial t} = \begin{array}{l} \sup_\varphi \left\{ \int_\Omega \{\varphi(z)\}^q \left\{\Phi^{(i)}(t,z) - \Phi^{(i)}(t,x)\right\}\mathrm{d}z - \int_\Omega \phi(\varphi(z))\mathrm{d}z \right\} \\ -\delta\Phi^{(i)}(t,x) + wU_i\left(x,\{\mu(t,\cdot)\}\right) \end{array}, \quad t\in(0,T),\quad x\in\Omega \quad (28)$$

subjected to the terminal condition $\Phi^{(i)}(T,x)=0$. Herein, $\{\mu(t,\cdot)\} = \{\mu_i(t,\cdot)\}_{i\in\Xi}$ is the law of $\left\{X_t^{(i)}\right\}_{i\in\Xi}$ and the supremum in (28) is considered with respect to all non-negative and measurable mapping $\varphi$ such that $\int_\Omega \varphi(z)\mathrm{d}z = 1$. If a maximizing $\varphi$ in (28) exists at $(t,x)$, it is denoted as $\varphi_i^*(t,x,z)$. Then, the time-forward PIDE, the Fokker–Planck (FP) equation associated with the MFG, for $i\in\Xi$ is given by

$$\frac{\mathrm{d}}{\mathrm{d}t}\mu_i(t,A) = \int_{x\in A}\int_{y\in\Omega} \left\{\varphi_i^*(t,y,x)\right\}^q \mu_i(t,\mathrm{d}y)\mathrm{d}x - \int_{x\in A}\int_{y\in\Omega} \left\{\varphi_i^*(t,x,y)\right\}^q \mathrm{d}y\,\mu_i(t,\mathrm{d}x) \quad \text{for any } A\in\mathcal{B} \quad (29)$$

subjected to an initial condition $\mu_i(0)\in\mathcal{P}_{m_i}$, which resembles the PCD (16) when $\rho_i = \left(\varphi_i^*\right)^q$. The mutually coupled PIDEs (28) and (29) require solutions in opposite directions in time with each other.

### 3.2 Specific case for logit dynamics

The coefficient $\phi$ is explicitly specified here. The strategy followed for this study is based on the Tsallis formalism of maximization problems subjected to entropic penalization [54,60,61], which designs $\phi$ as a generalized divergence:

$$\phi(u) = \eta\times\begin{cases} \frac{1}{1-q}\left(1-u^q+q(u-1)\right) & (q\neq 1) \\ u\ln u - u + 1 & (q=1) \end{cases} \text{ with } \frac{\mathrm{d}\phi(u)}{\mathrm{d}u} = \eta\times\begin{cases} \frac{q}{1-q}\left(1-u^{q-1}\right) & (q\neq 1) \\ \ln u & (q=1) \end{cases},\quad u\geq 0 \quad (30)$$

with the convention $0\ln 0 = 0$. The function $\phi=\phi(u)$ is non-negative, smooth, convex, and has the global minimum value 0 at $u=1$. The minimum corresponds to the uniform distribution case $\varphi(\cdot)=1$, with which players randomly update their actions without any reference to the current action. Hence, no control cost is incurred. Some positive cost can be incurred when $\varphi$ become positive. The case $q=1$ corresponds to the case where $\phi$ is proportional to the classical relative entropy [55], and the case $q > 1$ ($q \in (0,1)$) more strongly (more weakly) penalizes the control value deviating from 1.

The maximization problem in (28) is formally solved as

$$\begin{aligned} \varphi_i^*(t,x,z) &= \arg\max_{\varphi} \left\{ \int_\Omega \{\varphi(z)\}^q \left\{\Phi^{(i)}(t,z) - \Phi^{(i)}(t,x)\right\} \mathrm{d}z - \int_\Omega \phi(\varphi(z)) \mathrm{d}z \right\} \\ &= \frac{\exp_q\left(\eta^{-1}\Delta\Phi^{(i)}(t,z,x)\right)}{\int_\Omega \exp_q\left(\eta^{-1}\Delta\Phi^{(i)}(t,y,x)\right)\mathrm{d}y} \end{aligned} \quad (31)$$

at each $t \in (0,T)$ and $x \in \Omega$, where $\Delta\Phi^{(i)}(t,y,x) = \Phi^{(i)}(t,y) - \Phi^{(i)}(t,x)$. The FP equation (29) is identified as the GLD (6) when $\Delta U_i$ is replaced by $\Delta\Phi_i$. By (31), the HJB equation (28) becomes

$$\begin{aligned} -\frac{\partial\Phi^{(i)}(t,x)}{\partial t} &= \eta \ln_q\left(\int_\Omega \exp_q\left(\eta^{-1}\Delta\Phi^{(i)}(t,z,x)\right)\mathrm{d}z\right) \\ &\quad -\delta\Phi^{(i)}(t,x) + wU_i\left(x,\{\mu(t,\cdot)\}\right) \end{aligned}, \quad t \in (0,T) \text{ and } x \in \Omega. \quad (32)$$

As aforementioned, a similarity exists between the FP equation (29) and the GLD (6). Arguing that if $w = \delta$, then the MFG formally reduces to the GLD under the large-discount limit $\delta \to +\infty$ at which players are myopic and decide their actions based on their current state. Indeed, the equality (32) is rewritten as

$$\Phi^{(i)}(t,x) - U_i\left(x,\{\mu(t,\cdot)\}\right) = \frac{\eta \ln_q\left(\int_\Omega \exp_q\left(\eta^{-1}\Delta\Phi^{(i)}(t,z,x)\right)\mathrm{d}z\right) + \frac{\partial\Phi^{(i)}(t,x)}{\partial t}}{\delta}. \quad (33)$$

If the absolute value of the numerator on the right-hand side of (33) is bounded irrespective to $\delta > 0$, then the convergence $\Phi^{(i)} \to U_i$ is obtained, and the FP equation (29) reduces to the GLD (6). This argument is heuristic because the boundedness of each term in the numerator on the right-hand side of (33) has to be checked along with the convergence analysis of $\{\mu\}$ under the same limit. The issue is computationally addressed later. An obstacle in computation can be the finiteness of the terminal time $T$, with which the convergence $\Phi^{(i)} \to U_i$ may not hold well for small and large $t$ due to the influence of initial and terminal conditions of $\{\mu\}$ and $\Phi^{(i)}$, respectively. Nevertheless, for a sufficiently large $\delta > 0$, the convergence is expected to be approximately satisfied at an intermediate time $0 << t << T$ considering the turnpike property [62,63]. A stationary state would be approximately achieved under such conditions.

Hereafter, $w = \delta$ is assumed to consider the vanishing-discount limit.

***Remark 2*** The PIDEs (29) and (32) are satisfied both inside and on the boundary of the domain except at the initial and terminal times, respectively. Therefore, the PIDEs themselves serve as boundary conditions.

### 3.3 Theoretical analysis

Full analysis of the MFG model is beyond the scope of this study; however, the well-posedness of the HJB equation (28) is considered for any prescribed $\mu_i(t) \in \mathcal{P}_{m_i}$ ($t \ge 0$, $i \in \Xi$).

A collection of all bounded measurable functions on $[0,T]\times\Omega$ is denoted as $B([0,T]\times\Omega)$. First, the well-posedness of the following modified HJB equation from a viscosity viewpoint is studied:

$$H_i\left(t, x, \Phi^{(i)}(t,x), \Phi^{(i)}(t,\cdot), \frac{\partial \Phi^{(i)}(t,x)}{\partial t}\right) = 0 \tag{34}$$

with the Hamiltonian $H_i : [0,T]\times\Omega\times\mathbb{R}\times B([0,T]\times\Omega)\times\mathbb{R}\to\mathbb{R}$ given for each $i \in \Xi$ as

$$\begin{aligned} H_i(t,x,u,f(t,\cdot),Q) &= -Q + \delta u - \delta U_i(x,\{\mu(t,\cdot)\}) \\ &\quad -\eta \ln_q\left(\int_\Omega \exp_q\left(\eta^{-1}\min\{2L, \max\{-2L, f(t,z)-u\}\}\right)\mathrm{d}z\right) \end{aligned}, \tag{35}$$

where the minmax portion in (35) is introduced to well-define the non-local part considering **Assumption 2**, which turns out to be superficial for proper viscosity solutions to (35) (**Proposition 2**). Herein, $\{\mu(t,\cdot)\}$ is assumed such that $U_i(\cdot,\{\mu(t,\cdot)\})$ as a function of $0 \le t \le T$ is continuous.

A collection of all functions, continuous, lower-semicontinuous, upper-semicontinuous, and smooth on $[0,T]\times\Omega$ are denoted as $C([0,T]\times\Omega)$, $LSC([0,T]\times\Omega)$, $USC([0,T]\times\Omega)$, and $C^1([0,T]\times\Omega)$, respectively. Viscosity solutions to the modified HJB equation (34) are defined as follows.

***Definition 1***

***Viscosity super-solution***

*A collection* $\{\bar{\Phi}^{(i)}\}_{i\in\Xi}$ *of functions* $\bar{\Phi}^{(i)} \in LSC([0,T]\times\Omega)$ *is said to be a viscosity super-solution if it satisfies the following condition: for any test function* $\Psi^{(i)} \in C^1([0,T]\times\Omega)$ *(*$i = 1,2,3,\ldots,I$*) such that* $\bar{\Phi}^{(i)} - \Psi^{(i)}$ *is minimized at* $(i,t,x) = (\hat{i},\hat{t},\hat{x}) \in \Xi\times[0,T]\times\Omega$, *with* $\bar{\Phi}^{(\hat{i})}(\hat{t},\hat{x}) = \Psi^{(\hat{i})}(\hat{t},\hat{x})$, *it follows that*

$$H_{\hat{i}}\left(\hat{t}, \hat{x}, \bar{\Phi}^{(i)}(\hat{t},\hat{x}), \bar{\Phi}^{(i)}(\hat{t},\cdot), \frac{\partial \Psi^{(i)}(\hat{t},\hat{x})}{\partial t}\right) \ge 0 \quad \text{if } \hat{t} < T \tag{36}$$

*and* $\bar{\Phi}^{(i)}(T,\hat{x}) \ge 0$ *if* $\hat{t} = T$.

***Viscosity sub-solution***

*A collection* $\{\underline{\Phi}^{(i)}\}_{i\in\Xi}$ *of functions* $\underline{\Phi}^{(i)} \in USC([0,T]\times\Omega)$ *is said to be a viscosity sub-solution if it satisfies the following condition: for any test function* $\Psi^{(i)} \in C^1([0,T]\times\Omega)$ *(*$i = 1,2,3,\ldots,I$*) such that* $\underline{\Phi}^{(i)} - \Psi^{(i)}$ *is maximized at* $(i,t,x) = (\hat{i},\hat{t},\hat{x}) \in \Xi\times[0,T]\times\Omega$, *with* $\underline{\Phi}^{(\hat{i})}(\hat{t},\hat{x}) = \Psi^{(\hat{i})}(\hat{t},\hat{x})$, *it follows that*

$$H_i\left(\hat{t},\hat{x},\underline{\Phi}^{(i)}\left(\hat{t},\hat{x}\right),\underline{\Phi}^{(i)}\left(\hat{t},\cdot\right),\frac{\partial\Psi^{(i)}\left(\hat{t},\hat{x}\right)}{\partial t}\right)\leq 0 \quad \text{if} \quad \hat{t}<T \tag{37}$$

*and* $\underline{\Phi}^{(i)}\left(T,\hat{x}\right)\leq 0$ *if* $\hat{t}=T$.

***Viscosity solution***

*A collection* $\left\{\Phi^{(i)}\right\}_{i\in\Xi}$ *of functions* $\Phi^{(i)}\in C\left(\left[0,T\right]\times\Omega\right)$ *is said to be a viscosity solution if it is a viscosity super-solution and a viscosity sub-solution.*

The property of the Hamiltonian $H_i$ for each $i\in\Xi$ is found as follows:

✓ $H_i$ is Lipschitz continuous with respect to the last argument, increasing with respect to the third argument, such that for each $\left(t,x,f,Q\right)\in\left[0,T\right]\times\Omega\times B\left(\left[0,T\right]\times\Omega\right)\times\mathbb{R}$, it follows that

$$H_i\left(t,x,u_1,f\left(t,\cdot\right),Q\right)\leq H_i\left(t,x,u_2,f\left(t,\cdot\right),Q\right) \quad \text{if} \quad u_1\leq u_2 \quad \text{and the equality if and only if} \quad u_1=u_2, \tag{38}$$

✓ $H_i$ decreases with respect to the fourth argument, such that for each $\left(t,x,u,Q\right)\in\left[0,T\right]\times\Omega\times\mathbb{R}\times\mathbb{R}$, it follows that

$$H_i\left(t,x,u,f_1\left(t,\cdot\right),Q\right)\geq H_i\left(t,x,u,f_2\left(t,\cdot\right),Q\right) \quad \text{if} \quad f_1\leq f_2 \quad \text{on} \quad \left[0,T\right]\times\Omega. \tag{39}$$

The findings combined with the homogeneous terminal condition $\Phi^{(i)}=0$ imply that the modified HJB equation (34) admits at most one viscosity solution. Moreover, if a solution exists, it is bounded in $\left[-L,L\right]$.

***Proposition 2***

*Assume that* $U_i\left(x,\left\{\mu\left(t,\cdot\right)\right\}\right)$ *is uniformly continuous at all* $\left(t,x\right)\in\left[0,T\right]\times\Omega$. *For any viscosity super-solution* $\left\{\overline{\Phi}^{(i)}\right\}_{i\in\Xi}$ *and viscosity sub-solution* $\left\{\underline{\Phi}^{(i)}\right\}_{i\in\Xi}$, *it follows that* $\overline{\Phi}^{(i)}\geq\underline{\Phi}^{(i)}$ *on* $\left[0,T\right)\times\Omega$ *for any* $i\in\Xi$. *Moreover, any viscosity solution* $\left\{\Phi^{(i)}\right\}_{i\in\Xi}$ *satisfies* $-L\leq\Phi^{(i)}\leq L$ *on* $\left[0,T\right]\times\Omega$ *for any* $i\in\Xi$.

Bounds including those in **Proposition 2** are satisfied by numerical solutions to the discretization methods of this study as shown in the next section.

## 4. Computational analysis

### 4.1 Numerical scheme

A finite difference scheme to discretize the coupled PIDEs (29) and (32) in the MFG is discussed. The discretization employs a naïve way (e.g., Chapter 5 in [51]) that utility, probability measure, and value function are discretized at each discrete points distributed uniformly over the domain. The discretized equations are solved recursively as shown in **Algorithm 1**. The numerical method for the GLD is simpler as explained in **Remark 4**. The minmax operation as in (35) is not applied to the discretization as it turns out to be unnecessary (**Proposition 3**).

First, the space-time domain $[0,T]\times\Omega$ is discretized by vertices $\mathrm{P}_{k,l}:(t_k,x_l)$ with $t_k=k\Delta t$ ($k=0,1,2,...,N_t$) and $x_l=(l-1/2)\Delta x$ ($l=1,2,3,...,N_x$), where $\Delta t=T/N_t$ and $\Delta x=1/N_x$ for some $N_t,N_x\in\mathbb{N}$. The cell $\mathrm{C}_l$ associated with the vertex $\mathrm{P}_{\cdot,l}$ is given as $\mathrm{C}_l=[(l-1)\Delta x,l\Delta x)$ ($l=1,2,3,...,N_x-1$) and $\mathrm{C}_{N_x}=[(N_x-1)\Delta x,1]$. The discretized $\Phi^{(i)}$ at $\mathrm{P}_{k,l}$ is denoted by $\Phi^{(i)}_{k,l}$, and that on $\{t_k\}\times\mathrm{C}_l$ by $\mu_{i,k,l}$. Assume that the initial condition of discretized $\{\mu(0)\}$ is prepared for all $\{t_0\}\times\mathrm{C}_l$ ($l=1,2,3,...,N_x$) with $\sum_{l=1}^{N_x}\mu_{i,0,l}=m_i$ for each $i\in\Xi$, and the terminal condition of $\Phi^{(i)}$ is implemented as $\Phi^{(i)}_{N_t,\cdot}=0$.

Initially, the utility $U_i$ is discretized by assuming the generic form:

$$U_i\left(x,\{\mu(t,\cdot)\}\right)=F_i\left(x,\left\{\int_\Omega G_j(x)\mu_j(t,\mathrm{d}x)\right\}_{j=1,2,3...,I}\right),\quad x\in\Omega \tag{40}$$

with each $F_i:\Omega\times\mathbb{R}^I\to\mathbb{R}$ and $G_j:\Omega\to\mathbb{R}$ being continuous. Herein, (40) is discretized as

$$U_i\left(x,\{\mu(t,\cdot)\}\right)\approx U_{i,k,l}=F_i\left(x_l,\left\{\sum_{m=1}^{I}G_j(x_m)\mu_{j,k,m}\right\}_{j=1,2,3...,I}\right)\ \text{at}\ \mathrm{P}_{k,l}. \tag{41}$$

Based on (41), for $k=0,1,2,...,N_t-1$ and $l=1,2,3,...,N_x$, (32) is discretized as

$$-\frac{\Phi^{(i)}_{k+1,l}-\Phi^{(i)}_{k,l}}{\Delta t}=\eta\ln_q\left(\sum_{m=1}^{I}\exp_q\left(\eta^{-1}\left(\Phi^{(i)}_{k+1,m}-\Phi^{(i)}_{k+1,l}\right)\right)\Delta x\right)-\delta\Phi^{(i)}_{k+1,l}+\delta U_{i,k,l}, \tag{42}$$

which can be explicitly solved for $\Phi^{(i)}_{k,l}$ as

$$\Phi^{(i)}_{k,l}=\Phi^{(i)}_{k+1,l}+\Delta t\times\left\{\eta\ln_q\left(\sum_{m=1}^{I}\exp_q\left(\eta^{-1}\left(\Phi^{(i)}_{k+1,m}-\Phi^{(i)}_{k+1,l}\right)\right)\Delta x\right)-\delta\Phi^{(i)}_{k+1,l}+\delta U_{i,k,l}\right\}. \tag{43}$$

Then, $\varphi^*_i$ is discretized as

$$\varphi^*_i(t,z,x)\approx\varphi^*_{i,k,m,l}=\frac{\exp_q\left(\eta^{-1}\left(\Phi^{(i)}_{k,m}-\Phi^{(i)}_{k,l}\right)\right)}{\sum_{o=1}^{I}\exp_q\left(\eta^{-1}\left(\Phi^{(i)}_{k,o}-\Phi^{(i)}_{k,l}\right)\right)\Delta x}\quad\text{for}\ (t,x,z)=(t_k,x_m,x_l). \tag{44}$$

The FP equation (29) is discretized using (44) as

$$\frac{\mu_{i,k+1,l}-\mu_{i,k,l}}{\Delta t}=\sum_{m=1}^{N_x}\left(\varphi^*_{i,k,l,m}\right)^q\mu_{i,k,m}\Delta x-\left(\sum_{m=1}^{N_x}\left(\varphi^*_{i,k,m,l}\right)^q\Delta x\right)\mu_{i,k,l},\quad k=0,1,2,...,N_t-1,\ l=1,2,3,...,N_x, \tag{45}$$

which can be explicitly solved for $\mu_{i,k+1,l}$ as

$$\mu_{i,k+1,l}=\mu_{i,k,l}+\Delta t\left\{\sum_{m=1}^{N_x}\left(\varphi^*_{i,k,l,m}\right)^q\mu_{i,k,m}\Delta x-\left(\sum_{m=1}^{N_x}\left(\varphi^*_{i,k,m,l}\right)^q\Delta x\right)\mu_{i,k,l}\right\}. \tag{46}$$

The discretized MFG is solved to find $\Phi^{(i)}_{k,l}$ and $\mu_{i,k,l}$ for all $k = 0,1,2,...,N_t$, $l = 1,2,3,...,N_x$, and $i \in \Xi$. To solve (43) and (46), an iteration method is required as the equations are non-linearly coupled and must be integrated in the opposite directions of $k$. The issue is addressed by using **Algorithm 1**. Herein, subscript $[\cdot]$ represents iteration count and $\varsigma \in (0,1]$ is a relaxation constant. We use an abbreviation like $\{\mu_{i,k,l}\}_{i\in\Xi,k=0,1,2,...,N_t,l=1,2,3,...,N_x} = \{\mu_{i,k,l}\}$.

***Algorithm 1***

**0.** Prepare an initial guess $\{\mu_{i,k,l,[0]}\}$ and $\{\Phi^{(i)}_{k,l,[0]}\}$, and set $r = 1$ (initialization of the iteration count).

**1.** Compute $\{\Phi^{(i)}_{k,l,[r]}\}$ using (43) where $\mu_{i,k,l} = \mu_{i,k,l,[r-1]}$.

**2.** Compute $\{\mu_{i,k,l,[r]}\}$ using (46) where $\Phi^{(i)}_{k,l} = \Phi^{(i)}_{k,l,[r]}$.

**3.** Compute the iteration error $\varepsilon_r = \max\limits_{i,k,l}\left\{\left|\Phi^{(i)}_{k,l,[r]} - \Phi^{(i)}_{k,l,[r-1]}\right|, \left|\mu_{i,k,l,[r]} - \mu_{i,k,l,[r-1]}\right|\right\}$.

**4.** If $\varepsilon_r \le \bar{\varepsilon}\left(=10^{-10} \text{ in this study}\right)$, then obtain $\{\Phi^{(i)}_{k,l,[r]}\}$, $\{\mu_{i,k,l,[r]}\}$, the numerical solution, and terminate the algorithm.

**5.** Set $\Phi^{(i)}_{k,l,[r]} \to \varsigma\Phi^{(i)}_{k,l,[r]} + (1-\varsigma)\Phi^{(i)}_{k,l,[r-1]}$ and $\mu_{i,k,l,[r]} \to \varsigma\mu_{i,k,l,[r]} + (1-\varsigma)\mu_{i,k,l,[r-1]}$. Then, go to the step 1 with increasing $r \to r+1$.

**Algorithm 1** is a version of the alternating sweeping method with relaxation [64] that solves forward and backward equations alternatively. Selecting a smaller value of the parameter $\varsigma$ enhances iteration stability but possibly slows down its convergence. Theoretically, the convergence of the algorithm is achieved if $\varsigma$ is sufficiently small and the maximum eigen value of the system, as a collection of all discretized equations (42) and (46), is smaller than 1. However, verifying the condition is difficult in general. In this study, $\alpha = 0.5$ is sufficient to obtain stable numerical solutions.

***Remark 3*** The iteration needed to complete Algorithm was at most few tens and the convergence speed of the error was exponential in $r$. For example, we have $\varepsilon_r = -0.272r + 0.180$ ($R^2 = 0.999$) for second application in **Section 4**.

***Remark 4*** Discretization of the GLD considered in this study is

$$\mu_{i,k+1,l} = \mu_{i,k,l} + \Delta t\left\{\sum_{m=1}^{N_x}\left(\hat{\varphi}_{i,k,l,m}\right)^q \mu_{i,k,m}\Delta x - \left(\sum_{m=1}^{N_x}\left(\hat{\varphi}_{i,k,m,l}\right)^q \Delta x\right)\mu_{i,k,l}\right\} \tag{47}$$

with

$$\hat{\varphi}_{i,k,m,l} = \frac{\exp_q\left(\eta^{-1}\left(U_{i,k,m} - U_{i,k,l}\right)\right)}{\sum_{o=1}^{I}\exp_q\left(\eta^{-1}\left(U_{i,k,o} - U_{i,k,l}\right)\right)\Delta x} . \tag{48}$$

***Remark 5*** Existence of a solution to a semi-discretized version of the MFG is briefly discussed. For each $l = 1,2,3,...,N_x$, we consider the system of forward-backward ordinary differential equations

$$-\frac{\mathrm{d}\Phi_l^{(i)}(t)}{\mathrm{d}t} = \eta \ln_q\left(\sum_{m=1}^{I}\exp_q\left(\eta^{-1}\left(\Phi_m^{(i)}(t) - \Phi_l^{(i)}(t)\right)\right)\Delta x\right) - \delta\Phi_l^{(i)}(t) + \delta U_{i,l} , \quad 0 \le t < T , \tag{49}$$

$$\frac{\mathrm{d}\mu_i(t,C_l)}{\mathrm{d}t} = \sum_{m=1}^{N_x}\varphi^*_{i,l,m}\mu_i(t,C_m)\Delta x - \left(\sum_{m=1}^{N_x}\left(\varphi^*_{i,m,l}\right)^q \Delta x\right)\mu_i(t,C_l) , \quad 0 < t \le T \tag{50}$$

with suitable initial and terminal conditions, where

$$\varphi^*_{i,m,l} = \frac{\exp_q\left(\eta^{-1}\left(\Phi_m^{(i)}(t) - \Phi_l^{(i)}(t)\right)\right)}{\sum_{o=1}^{I}\exp_q\left(\eta^{-1}\left(\Phi_o^{(i)}(t) - \Phi_l^{(i)}(t)\right)\right)\Delta x} . \tag{51}$$

Existence of a solution to the system (49)-(50) follows due to the Lipschitz continuity of each $\varphi^*_{i,l,m}$ with respect to each $\mu_i$ by **Assumptions 1** and the uniform boundedness $\left|\Phi_l^{(i)}\right| \le L$ that follows essentially the same way with the Proposition 4 in Gomes et al. [65]. Ensuring uniqueness of solutions of the semi-discretized system requires additional assumptions such as a strict concavity of each $U_i$ with respect to each $\mu_i$.

We present propositions that reveal several properties of the numerical solutions. Moreover, without any loss of generality, we assume nonnegativity of $U_i$ as well as that of its discretized version $U_{i,k,l}$. This is satisfied if $F_i$ in (40) is nonnegative.

***Proposition 3***

*Given a sequence of* $\{\mu_{i,k,l}\}$, *such that each* $\mu_{i,k,l}$ *is non-negative and* $\sum_{l=1}^{N_x}\mu_{i,k,l} = m_i$ *for each* $i \in \Xi$ *and* $k = 0,1,2,...,N_t$, *if* $\Delta t \le \frac{1}{1+\delta}$ *then it follows that*

$$0 \le \Phi_{k,l}^{(i)} \le L \quad \textit{for all} \quad i \in \Xi , \quad k = 0,1,2,...,N_t , \textit{ and } l = 1,2,3,...,N_x . \tag{52}$$

***Proposition 4***

*Given a sequence of* $\{\Phi_{k,l}^{(i)}\}$ *such that each* $\Phi_{k,l}^{(i)}$ *satisfies the bound, if* $\Delta t \le \left(\frac{\exp_q\left(-2\eta^{-1}L\right)}{\exp_q\left(2\eta^{-1}L\right)}\right)^q$, *then it follows that*

$$\mu_{i,k,l} \geq 0 \quad \textit{for all} \quad i \in \Xi, \quad k = 0,1,2,...,N_t, \textit{ and } \quad l = 1,2,3,...,N_x \tag{53}$$

*and*

$$\sum_{l=1}^{N_x} \mu_{i,k,l} = m_i \quad \textit{for each} \quad i \in \Xi \quad \textit{and} \quad k = 0,1,2,...,N_t . \tag{54}$$

***Remark 6*** According to **Propositions 3 and 4**, **Algorithm 1** always satisfies the inequalities (52), (53), and (55) if $0 < \Delta t \leq \min\left\{ \frac{1}{1+\delta}, \left( \frac{\exp_q\left(-2\eta^{-1}L\right)}{\exp_q\left(2\eta^{-1}L\right)} \right)^q \right\}$. The upper-bound of $\Delta t$ is independent from $\Delta x$.

***Remark 7*** The convergence theorem (Theorem 2.1 [66]) suggests that numerical solutions to (43), given a sequence of discretized probability measures of $\{\mu\}$, converge to a viscosity solution to the modified HJB equation (34) if it exists. By **Proposition 2**, the solution is its unique viscosity solution. Full analysis of convergence of numerical solutions to the MFG can be explored further as it crucially depends on the regularity of utility $U_i$, as discussed in **Remark 5**.

Unless otherwise specified, resolutions $N_x = 150$ ($\Delta x = 1/150$ (1/day)) and $N_t = 240 \times 150$ ($\Delta t = 1/150$ (day), $T = 240$ (day)) are considered below. Herein, $x$ can be considered as the arrival intensity (1/day) of players at a site in each problem.

### 4.2 Target problems

Two applications of the GLD and its MFG counterpart, involving recreational fishing by legal and illegal anglers and sustainable tourism, are considered for this study. In both problems, $x$ represents a scaled arrival intensity of players and $\Xi = \{1,2\}$. Both models are simple and prepared for analyzing potential game-cases related to applications and are not intended to real-world applications.

#### 4.2.1 Recreational fishing

Fishing in inland water bodies, such as rivers and lakes, has been a major recreational activity worldwide [67,68]. Generally, each angler must pay some license fee per day or per year, depending on regional fishery cooperatives authorizing local fishery resources. Studies report an increasing harvesting pressure from anglers to fishes [69], suggesting the requirement for regulation to limit fish catch. Nakamura [70] reported that an estimated 71.8 % of anglers in riverine inland fisheries of Japan have been paying license fee for several major fishes. Thus, fishing pressure from illegal anglers (those not paying license fee) requires suppression.

For this study, two population types: legal ($i = 1$) and illegal anglers ($i = 2$), are considered. The utility of a concave potential game-type is proposed as

$$U_i(x,\mu_1,\mu_2) = \underbrace{x^{\alpha}}_{\text{Fishing utility}} - \underbrace{\underbrace{\beta x}_{\text{Proportional to arrival intensity}} \left(1+\kappa \underbrace{\frac{m_1}{m_i}\Delta_{i,2}}_{\text{Penalty}}\right) \underbrace{\int_\Omega y\left(\mu_1(\mathrm{d}y)+\mu_2(\mathrm{d}y)\right)}_{\text{Cost due to harvesting pressure}}}_{\text{Fishing cost}}, \tag{56}$$

where $\alpha > 0$ is a parameter modulating the fishing utility, $\beta > 0$ is the proportional cost, $\kappa > 0$ is a penalty to be incurred for the illegal anglers, and $\Delta_{i,j}$ is the Kronecker Delta ($\Delta_{i,j} = 1$ if $i = j$ and $\Delta_{i,j} = 0$, otherwise). The second term in (56) denotes the net harvesting cost that accounts for the monetary cost, which is proportional to the arrival intensity through the license fee ($\beta x$), as some penalty ($\kappa$) can be incurred for illegal anglers when detected. The detection is more likely to occur if the ratio $\frac{m_1}{m_2}$ is large, and local fishery cooperative modulates the license fee to regulate fish catch by anglers through the average arrival rate $\int_\Omega y\left(\mu_1(\mathrm{d}y)+\mu_2(\mathrm{d}y)\right)$. Indeed, Van Leeuwen [71] reported that illegal incidents were negatively correlated with the number of licensed anglers. The potential function associated with the utility (56) is found as (**Proposition 5** [44])

$$\mathfrak{P}(\mu_1,\mu_2) = \int_\Omega y^{\alpha}\left(\mu_1(\mathrm{d}y)+\mu_2(\mathrm{d}y)\right) - \frac{\beta}{2}\left(1+\kappa\frac{m_1}{m_i}\Delta_{i,2}\right)\left\{\int_\Omega y\left(\mu_1(\mathrm{d}y)+\mu_2(\mathrm{d}y)\right)\right\}^2. \tag{57}$$

Notably, the Fréchet derivative is formally defined as $\frac{\delta \mathfrak{P}}{\delta \mu_i} = U_i$. Cheung and Lahkar [44] reported that the potential of the function (57) is concave, and the set of Nash equilibria of the utility (56) is convex and coincides with the set of maximizing pairs $(\mu_1,\mu_2)$ of $\mathfrak{P}$.

### 4.2.2 Sustainable tourism

As surveyed by Pokki et al. [72], anglers visiting a water body for recreational purpose can be categorized into two clusters; the first cluster for residents and the second for fishing tourists who travel longer and arrive less frequently. On ecotourism, Grilli et al. [73] suggested that tourists can be categorized into three clusters, including tourists who strongly concern sustainable tourism, those less concern sustainable tourism, and those indifferent to sustainable tourism. Additionally, heterogeneous preferences among tourists of ecotourism have been reported by Shi et al. [74], where the behavior of tourists was modeled based on a heterogeneous model. Antoci et al. [75] considered an evolutionary game for analyzing overtourism model containing two populations, residents, and tourists, and explored stable equilibria at which both populations coexisted. For example, a survey has been conducted in the Shiramine Village of the Hakusan Mountainous area, Hokuriku region, Japan [23], located in the middle of the biosphere reserve of United Nations Educational, Scientific and Cultural Organization. The village is famous for hot springs and mountainous inland fisheries. According to Mammadova et al. [76], the village is currently facing depopulation, accessibility problems (especially for the elderly), and difficulties with tourism. The village residents are considering the importance of preserving environment and traditional culture. Trips to biological reserves, such as Sharmaine Village, should be balanced among local environment, human lives, and tourism.

The two population types: residents ($i=1$) and tourists ($i=2$) are considered. A regularized indicator function $\mathbb{I}_{\varepsilon}(\cdot \le \hat{x})$ is set with $\hat{x} \in (0,1)$ and $\varepsilon > 0$ as $\mathbb{I}_{\varepsilon}(x \le \hat{x}) = (1+\tanh((\hat{x}-x)/\varepsilon))/2$ ($x \in \mathbb{R}$), whose limit under $\varepsilon \to +0$ is a discontinuous step function. The problem in the second application has a utility that is of the potential type but more irregular than (57):

$$U_i(x,\mu_1,\mu_2) = \underbrace{\underbrace{\frac{1}{\theta + \int_{\Omega}(1-\mathbb{I}_{\varepsilon}(x \le \hat{x}))(\mu_1(\mathrm{d}y)+\mu_2(\mathrm{d}y))}}_{\text{Utility decrease by overtourism}} \mathbb{I}_{\varepsilon}(x \le \hat{x}) - \underbrace{\gamma_i x}_{\text{Travel cost}}}_{\text{Net Utility}}, \quad (58)$$

where $\theta > 0$ is a constant, $\gamma_i > 0$ denotes the travel cost with $\gamma_1 < \gamma_2$, and $\hat{x} \in (0,1)$ is the threshold value above which overtourism can occur. The utility of each population decreases if the expectation $\int_{\Omega}(1-\mathbb{I}_{\varepsilon}(x \le \hat{x}))\mu_i(\mathrm{d}y)$ that measures overtourism in the mean is large. The constant $\theta$ in the denominator of the first term on the right hand-side of (58) improves regularity of the utility and prevents divisions by zero, but the utility does not satisfy **Assumption 1**. The utility (58) implies that the travel cost is larger for the tourists than residents, and utility of the biological reserve is possibly lost if there are any excessive populations. The potential function $\mathfrak{P}$ corresponding to (58) is given by

$$\mathfrak{P}(\mu_1,\mu_2) = -\ln\left(\theta + \int_{\Omega}(1-\mathbb{I}_{\varepsilon}(x \le \hat{x}))(\mu_1(\mathrm{d}y)+\mu_2(\mathrm{d}y))\right) - \int_{\Omega} y(\gamma_1\mu_1(\mathrm{d}y)+\gamma_2\mu_2(\mathrm{d}y)). \quad (59)$$

The indicator function $\mathbb{I}_{\varepsilon}(\cdot \le \hat{x})$ is neither convex nor concave; hence, the potential function is not concave. Several values of $\varepsilon$ are computationally evaluated below.

### 4.3 Computational results: recreational fishing

First, the potential game of legal and illegal anglers is analyzed. Herein, $m_1 = 0.7$ and $m_2 = 0.3$ are set based on Nakamura [70]. Additionally, $\alpha = 0.5$ (the utility is concave with respect to the arrival intensity), $\beta = 2$, $\kappa = 0.1$, $q = 0.8$, $\eta = 0.01$, and $\delta = 1$ (1/day) are set. For visualization, the probability density $p_{i,k,l} = (\Delta x)^{-1}\mu_{i,k,l}$ is plotted in the sequel rather than the probability measure $\mu_{i,k,l}$. Numerical solutions to the GLD are judged as stationary when the difference of computed probability densities between successive time steps becomes smaller than $10^{-10}$. The initial conditions of the probability densities are $m_i(1+0.25(2x-1))$ and uniform distributions $m_i$ in GLD and MFG, respectively. Different initial conditions are selected to computationally suggest that stationary probability densities between GLD and MFG coincide for different initial conditions. Herein, a “stationary” solution to the MFG is at $t = T/2 = 120$ (day) considering the turnpike property discussed further.

The computed probability densities for the GLD are shown in **Fig. 1**, while the computed probability densities and value functions for the MFG are shown in **Fig. 2**. **Tables 1 and 2** present the convergence of numerical solutions (probability densities for GLD and both probability densities and value functions for MFG) as the computational resolution increases. The errors are evaluated by considering two norms: the maximum norm $\max_{l=1,2,3,\ldots,N_x} \left| u_{\text{num},l} - u_{\text{ref},l} \right|$ (**Table 1**) and average norm $\frac{1}{N_x}\sum_{l=1}^{N_x} \left| u_{\text{num},l} - u_{\text{ref},l} \right|$ (**Table 2**). Herein, $u$ represents the probability density or value function, "ref" denotes reference solution that is a numerical solution with $N_x = 300$ and $N_t = 240 \times 300$, and "num" denotes numerical solutions with coarser resolution. The center points of cells do not coincide between reference and coarser numerical solutions, and hence the cell averages of the latter are computed in coarser cells to compute the error norms. According to **Figs. 1 and 2**, profiles of the computed probability densities are smooth for both $i = 1, 2$, and the convergence check with respect to the maximum and average norms work well.

The numerical solutions to the MFG exhibit turnpike, namely the probability densities are almost constant for intermediate values of time $t$, such as $t = 120$ (**Fig. 2**). Considering this finding, the computed "stationary" probability densities of the GLD are compared with those of the MFGs for different values of the discount rate $\delta$ (**Fig. 3**). The computational results visually demonstrate that the probability densities at the middle of the turnpike converge to those of GLD with an increase in $\delta$. **Table 3** presents the convergence of the former to the latter in the order of $O\left(\delta^{-1}\right)$; hence, the error between the probability densities of the GLD and MFG decreases as $\delta$ increases, supporting the heuristic argument about the convergence of the latter to the former. Accordingly, the GLD can be obtained as a myopic limit of the corresponding MFG as per the computational verification results.

Furthermore, we analyze influences of the population of illegal anglers increase. The computed stationary probability densities of the MFG for three scenarios (Legal decrease: $m_1 = 0.5$, Neutral: $m_1 = 0.7$, Legal increase: $m_1 = 0.9$) corresponding to different values of the population $m_1$ of legal anglers are shown in **Fig. 4**. The decrease in legal angler population advocates harvesting by illegal anglers due to the decrease in penalization by not paying license fee. The probability density profile of the legal angler population does not change qualitatively. The increase in legal angler population leads to the illegal angler population concentrated at $x = 0$; hence, harvesting by illegal anglers is effectively ruled out. Simultaneously, the maxima of the arrival intensity distribution of legal anglers slightly shifts toward a smaller $x$ because the increase in harvesting pressure from the populations reduces the net utility.

**Table 1.** Convergence of numerical solutions to GLD and MFG with $(N_t, N_x) = (120\times m, m)$ relative to the maximum norm.

| | GLD | | MFG | | | |
|---|---|---|---|---|---|---|
| $m$ | $p_1$ | $p_2$ | $p_1$ | $p_2$ | $\Phi^{(1)}$ | $\Phi^{(2)}$ |
| 50 | 2.88.E-03 | 4.00.E-01 | 1.26.E-03 | 1.30.E-01 | 2.64.E-03 | 3.20.E-03 |
| 100 | 5.30.E-04 | 8.54.E-02 | 2.68.E-04 | 8.00.E-02 | 1.57.E-03 | 1.89.E-03 |
| 150 | 1.96.E-04 | 3.95.E-02 | 9.92.E-05 | 4.21.E-02 | 9.92.E-04 | 1.19.E-03 |

**Table 2.** Convergence of numerical solutions to GLD and MFG with $(N_t, N_x) = (120\times m, m)$ relative to the average norm.

| | GLD | | MFG | | | |
|---|---|---|---|---|---|---|
| $m$ | $p_1$ | $p_2$ | $p_1$ | $p_2$ | $\Phi^{(1)}$ | $\Phi^{(2)}$ |
| 50 | 9.55.E-04 | 2.12.E-02 | 3.76.E-04 | 6.58.E-03 | 9.98.E-05 | 2.70.E-04 |
| 100 | 1.82.E-04 | 3.66.E-03 | 9.27.E-05 | 2.35.E-03 | 2.47.E-05 | 6.70.E-05 |
| 150 | 4.83.E-05 | 9.21.E-04 | 3.58.E-05 | 9.18.E-04 | 1.05.E-05 | 2.83.E-05 |

**Table 3.** Difference of numerical solutions between GLD and MFG for different values of $\delta$ relative to the maximum norm ($\eta = 0.001$).

| $\delta$ | $p_1$ | $p_2$ |
|---|---|---|
| 1 | 1.68.E+00 | 2.59.E+00 |
| 5 | 4.43.E-01 | 6.17.E-01 |
| 10 | 2.31.E-01 | 3.16.E-01 |
| 25 | 9.49.E-02 | 1.28.E-01 |
| 50 | 4.78.E-02 | 6.44.E-02 |
| 100 | 2.40.E-02 | 3.23.E-02 |

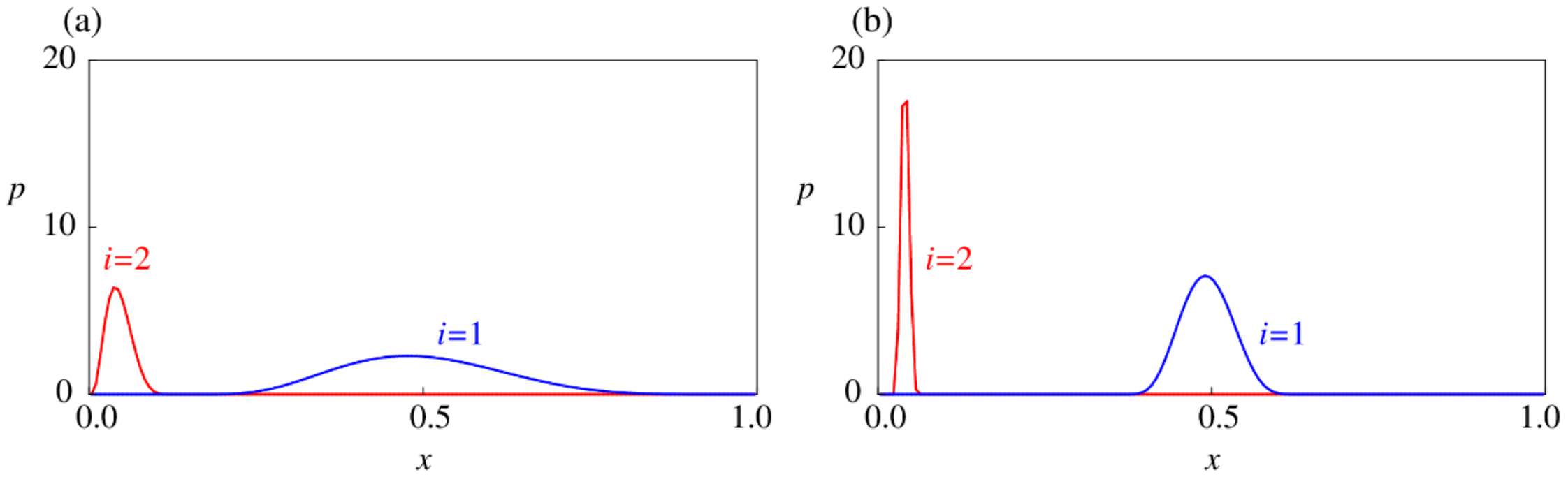

**Fig. 1.** Computed probability densities of GLD: (a) $\varepsilon = 0.01$ and (b) $\varepsilon = 0.001$.

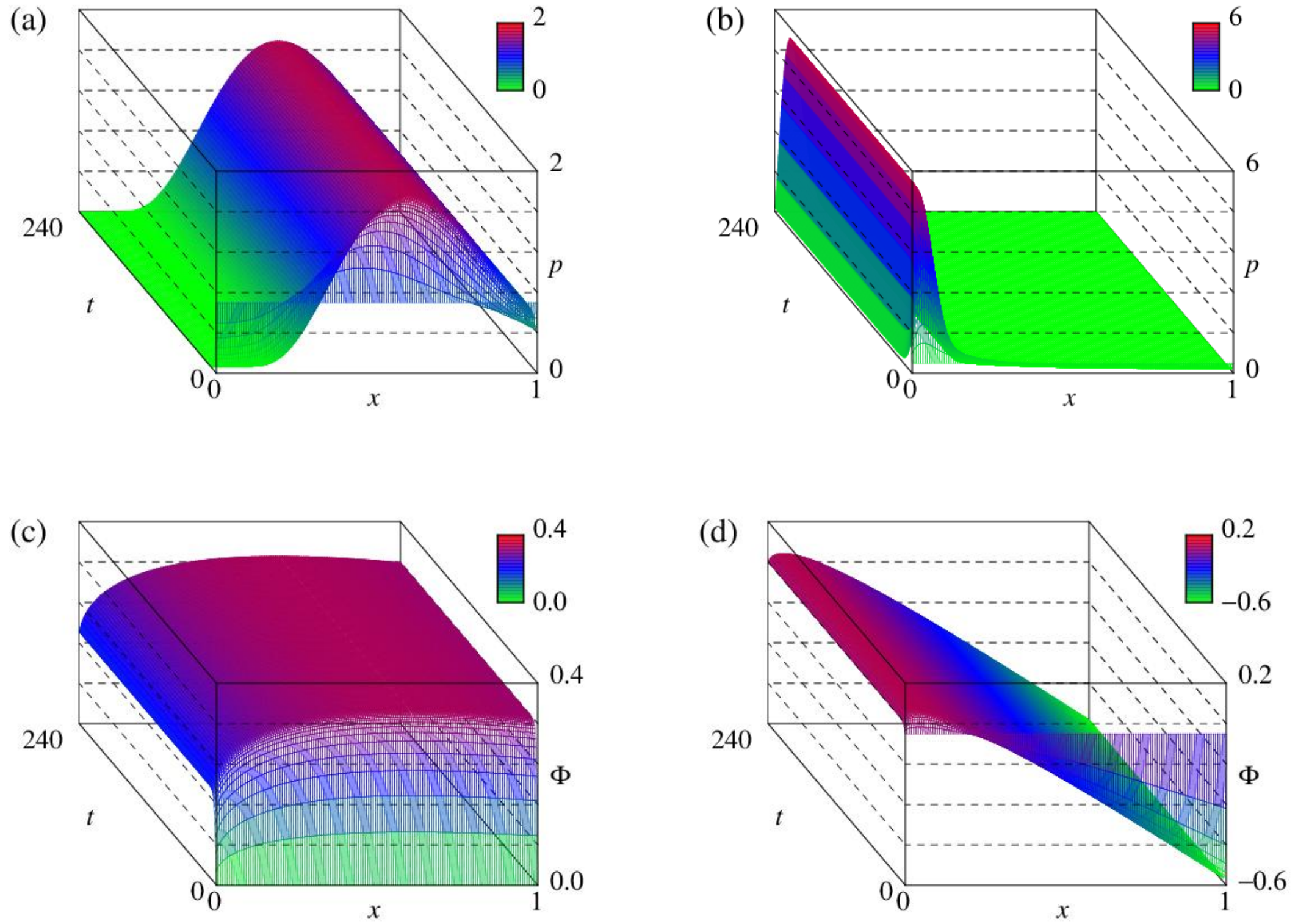

**Fig. 2.** Computed probability densities ((a) $i = 1$, (b) $i = 2$) and value functions ((c) $i = 1$, (d) $i = 2$) of MFG with $\delta = 1$.

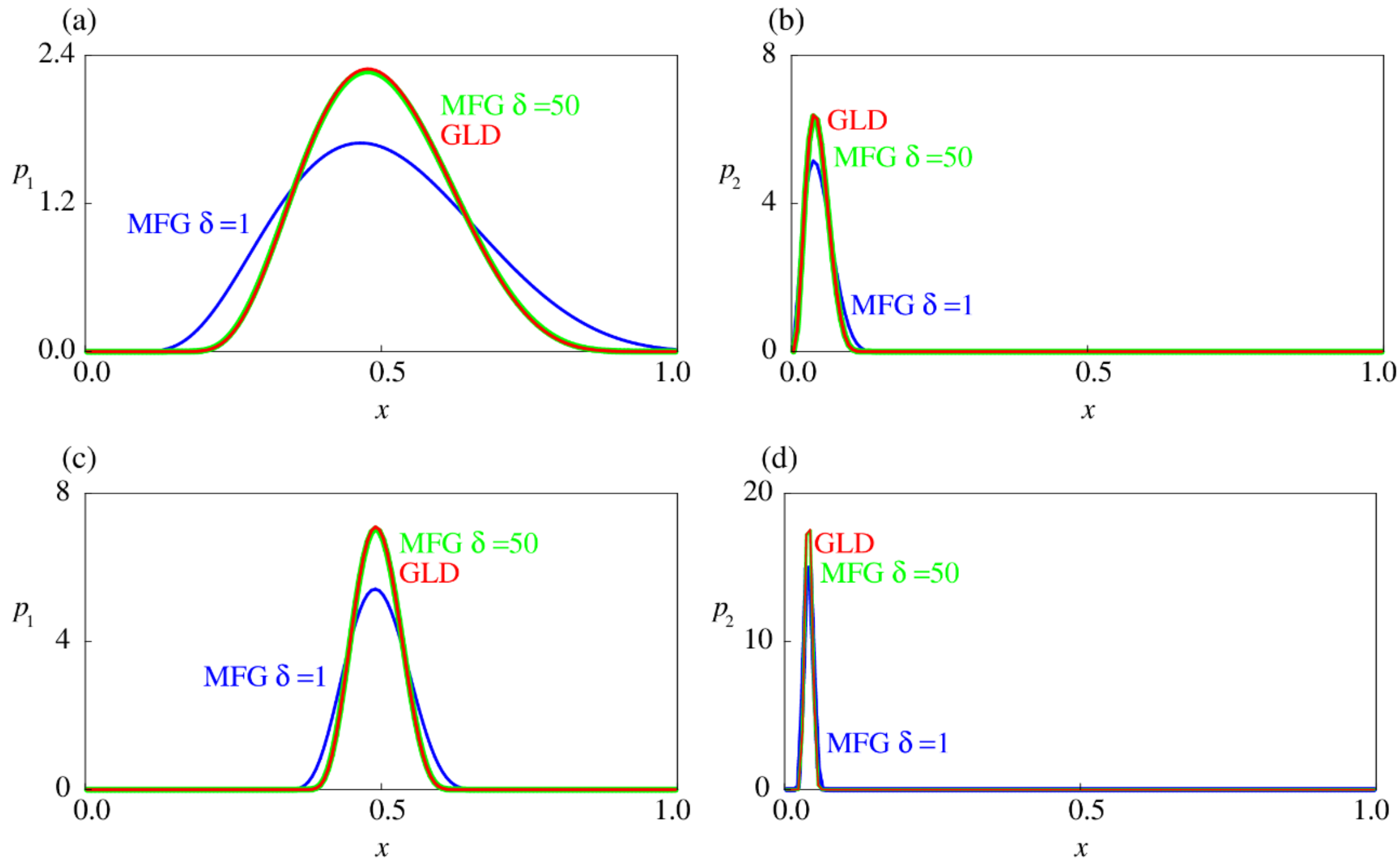

**Fig. 3.** Comparison of computed stationary probability densities between GLD and MFG with different values of $\delta$: $\eta = 0.01$ ((a) $i = 1$ and (b) $i = 2$) and $\eta = 0.001$ ((c) $i = 1$ and (d) $i = 2$).

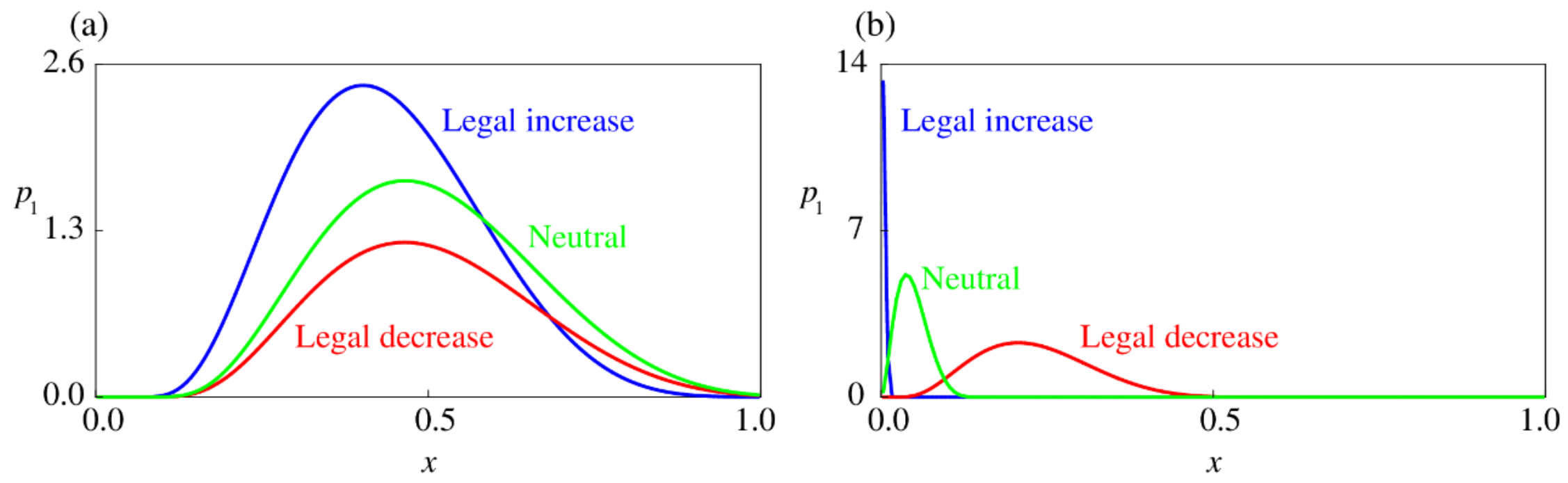

**Fig. 4.** Computed probability densities for GLD and MFG with different values of $m_1$ (Legal decrease: $m_1 = 0.5$, Neutral: $m_1 = 0.7$, Legal increase: $m_1 = 0.9$).

### 4.4 Computational results: Sustainable tourism

Unless other specified, the parameter values are: $\theta = 1$, $m_1 = 0.8$, and $m_2 = 0.2$ (tourist population is smaller than resident one); $\varepsilon = 10^{-6}$ (1/day), $\gamma_1 = 0.01$, and $\gamma_2 = 0.1$ (travel cost of tourists is larger than thar of residents); $\hat{x} = 0.65$, $q = 0.8$, $\eta = 0.01$, and $\delta = 1$ (1/day) are considered. The utility is continuous but with a large Lipschitz constant $L = O\left(\varepsilon^{-1}\right) = O\left(10^{6}\right)$; it is (almost) discontinuous.

The computed probability densities for the GLD are shown in **Fig. 5**, while the computed probability densities and value functions for the MFG are shown in **Fig. 6**. **Tables 3 and 4**, as in **Tables 1 and 2**, show the convergence of numerical solutions for the maximum and average norms. Profiles of the computed probability densities of the GLD are smooth for $i = 2$ but not for $i = 1$, suggesting that the convergence check with respect to the maximum norm does not work well due to the sharp transition of solutions. The average norm works better than the maximum one. Both the cases $i = 1$ and $i = 2$ have discontinuous probability densities, and the maximum norm does not work well (**Fig. 7**). A similar observation applies to the value functions (**Fig. 6**). The low regularity of the utility slows down convergence of numerical solutions. The discontinuity in the probability density $p_2$ is found to disappear as the travel cost $\gamma_2$ increases, while the density itself is still computationally 0 for $x > \hat{x}$. According to **Fig. 7**, computed probability densities of the MFG get closer to that of the GLD as the discount rate increases even for a utility that varies sharply with respect to $x$ (**Fig. 6**). The computational results imply that the MFG formally converges to the GLD as the discount rate increases for both continuous and (almost) discontinuous utility cases, which is considered partly due to the uniform boundedness of the value function irrespective to the discount rate, as stated in the previous section. More rigorously, the partial differential term of the value function in the right-hand side of (33) should be estimated as they may not be bounded. Nevertheless, the turnpike of the obtained numerical solutions suggests that this partial differential term is small for intermediate values of time $t$, supporting the heuristic argument of this study.

Additionally, a remark about the case where tourist population is smaller than resident one with $m_1 = 0.2$ and $m_2 = 0.8$ is provided. Although not presented in this study, for the GLD, the computed probability densities satisfy the relationships $\frac{p_1|_{(m_1,m_2)=(0.8,0.2)}}{p_1|_{(m_1,m_2)=(0.2,0.8)}} \approx 4$ and $\frac{p_2|_{(m_1,m_2)=(0.8,0.2)}}{p_2|_{(m_1,m_2)=(0.2,0.8)}} \approx \frac{1}{4}$, suggesting that the parameters $m_1, m_2$ only modulate the height of the probability densities but not their shapes. Therefore, action profiles to be selected by players at an equilibrium do not critically depend on the population ratios for the cases examined in this study.

Finally, the probability densities of MFG are analyzed with different values of $\varepsilon$ (**Fig. 8**). Profiles of the probability densities with $\varepsilon = 10^{-6}$ are discontinuous and continuous for $i = 1$ and $i = 2$, respectively, smeared out as $\varepsilon$ increases, and becomes similar at $\varepsilon = 1$. The obtained results imply that residents and tourists share a similar preference when the utility varies slowly with respect to the arrival intensity $x$. In all the cases, the computed probability densities are close to 0 for the intensity above the sustainability threshold $\hat{x}$. Although beyond the scope of this study, surveys to estimate arrival intensity profiles of residents and tourists can assist modelling the utilities, particularly their regularity conditions. Such an information can be helpful for designing price and environmental tax for sightseeing and analyzing their robustness.

**Table 4.** Convergence of numerical solutions to GLD and MFG with $(N_t, N_x) = (120 \times m, m)$ relative to the maximum norm.

| | GLD | | MFG | | | |
|---|---|---|---|---|---|---|
| $m$ | $p_1$ | $p_2$ | $p_1$ | $p_2$ | $\Phi^{(1)}$ | $\Phi^{(2)}$ |
| 50 | 4.67.E-01 | 3.20.E-04 | 5.36.E-01 | 1.34.E-02 | 8.16.E-04 | 8.45.E-04 |
| 100 | 1.00.E-05 | 8.00.E-05 | 1.00.E-05 | 2.50.E-05 | 0.00.E+00 | 1.00.E-06 |
| 150 | 4.66.E-01 | 9.15.E-05 | 5.35.E-01 | 1.38.E-02 | 2.71.E-04 | 5.91.E-04 |

**Table 5.** Convergence of numerical solutions to GLD and MFG with $(N_t, N_x) = (120 \times m, m)$ relative to the average norm.

| | GLD | | MFG | | | |
|---|---|---|---|---|---|---|
| $m$ | $p_1$ | $p_2$ | $p_1$ | $p_2$ | $\Phi^{(1)}$ | $\Phi^{(2)}$ |
| 50 | 1.87.E-02 | 3.31.E-05 | 2.14.E-02 | 1.84.E-04 | 3.64.E-04 | 2.63.E-04 |
| 100 | 6.10.E-07 | 7.33.E-06 | 2.00.E-07 | 3.26.E-06 | 0.00.E+00 | 6.20.E-07 |
| 150 | 4.66.E-03 | 2.40.E-06 | 5.35.E-03 | 5.64.E-04 | 1.21.E-04 | 8.64.E-05 |

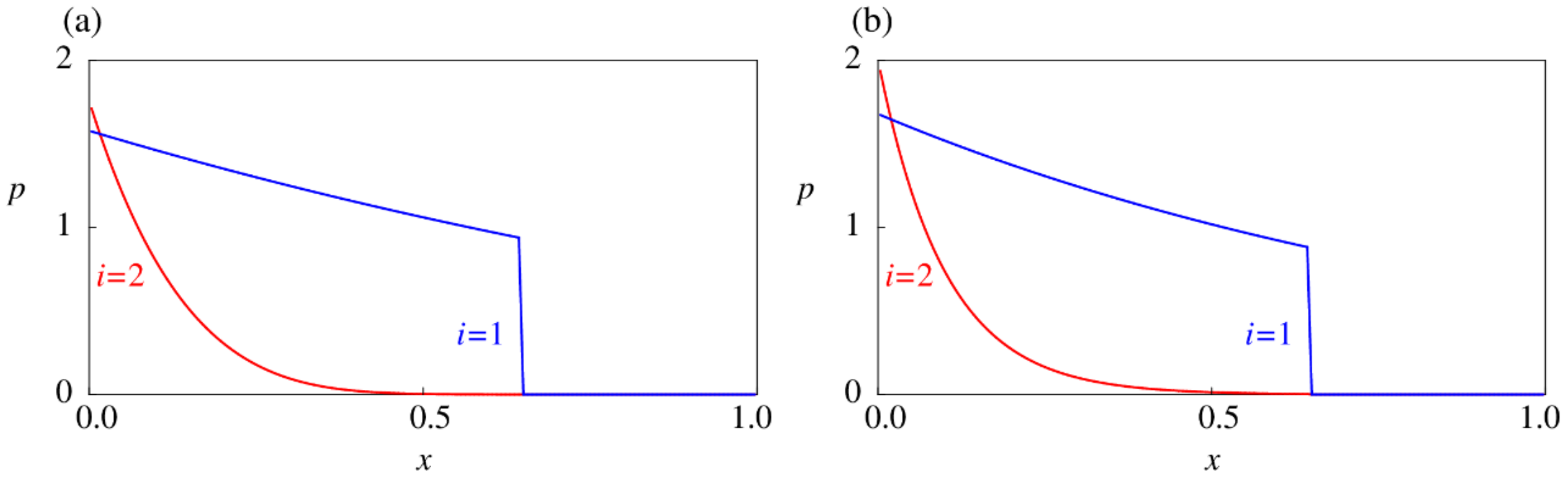

**Fig. 5.** Computed probability densities of GLD (a) $q = 0.8$ and (b) $q = 0.1$.

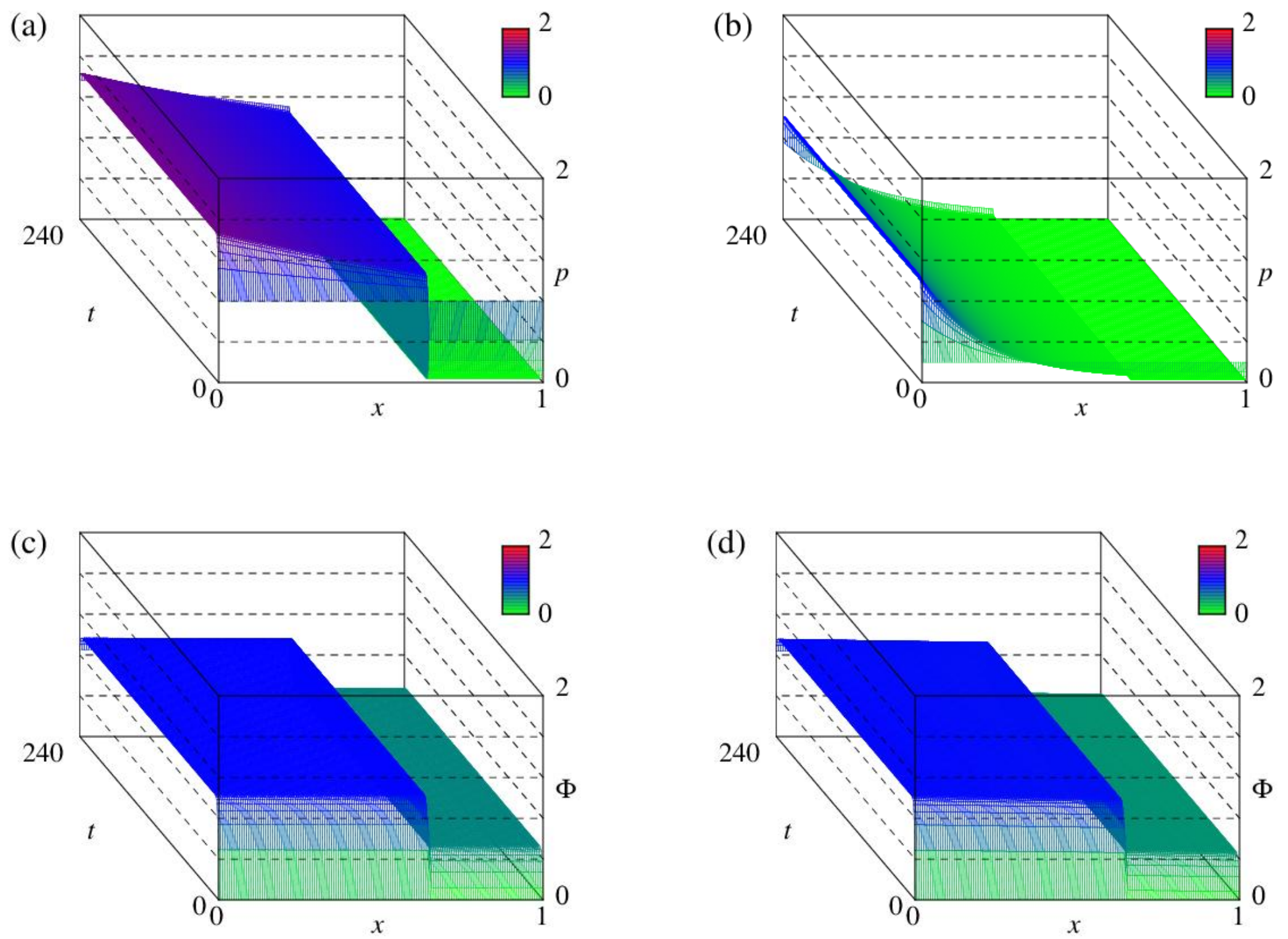

**Fig. 6.** Computed probability densities ((a) $i = 1$, (b) $i = 2$) and value functions ((c) $i = 1$, (d) $i = 2$) of MFG with $\delta = 1$.

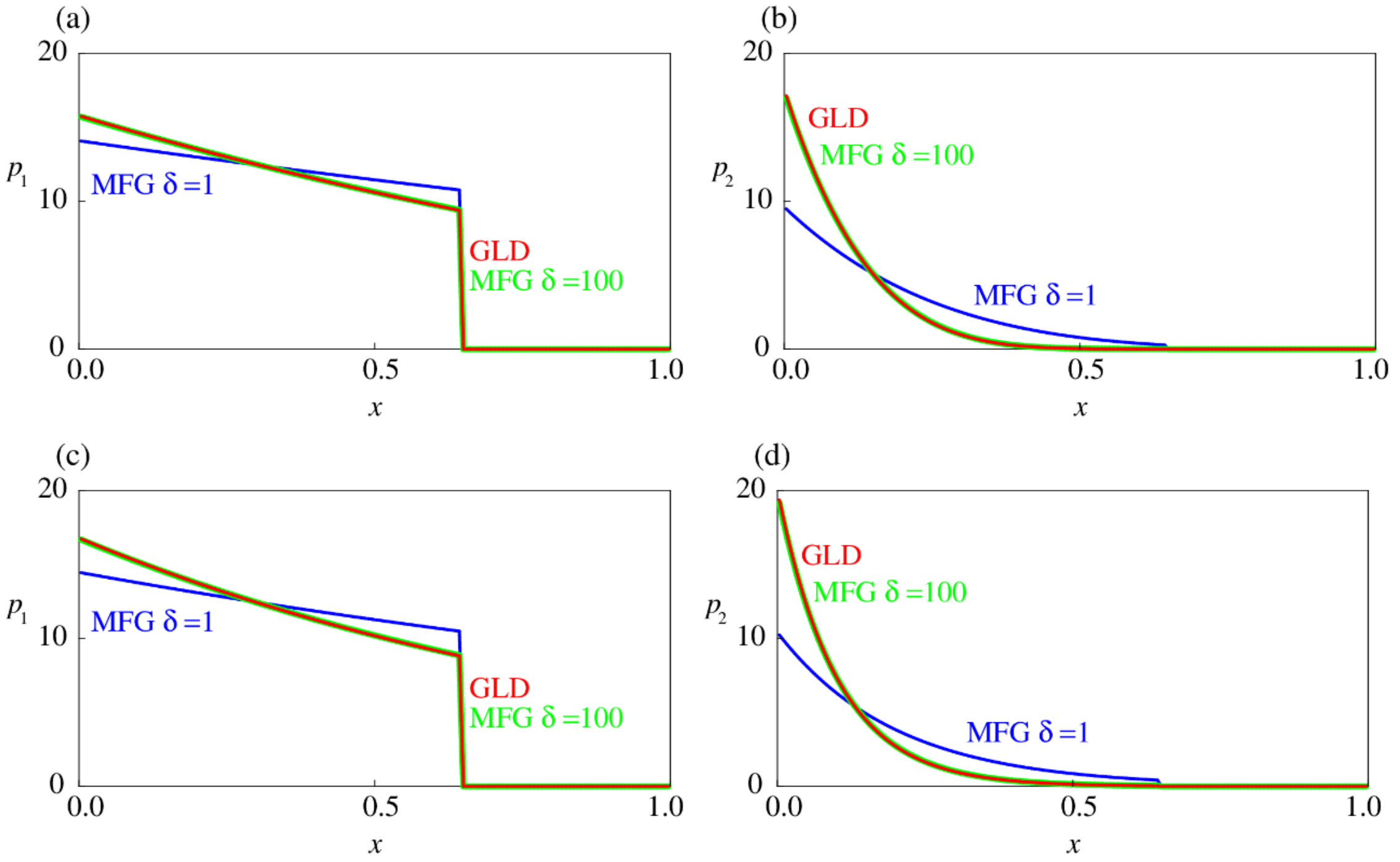

**Fig. 7.** Comparison of computed "stationary" probability densities between GLD and MFG with different values of $\delta$ : $q = 0.8$ ((a) $i = 1$ and (b) $i = 2$ ) and $q = 1$ ((c) $i = 1$ and (d) $i = 2$ ).

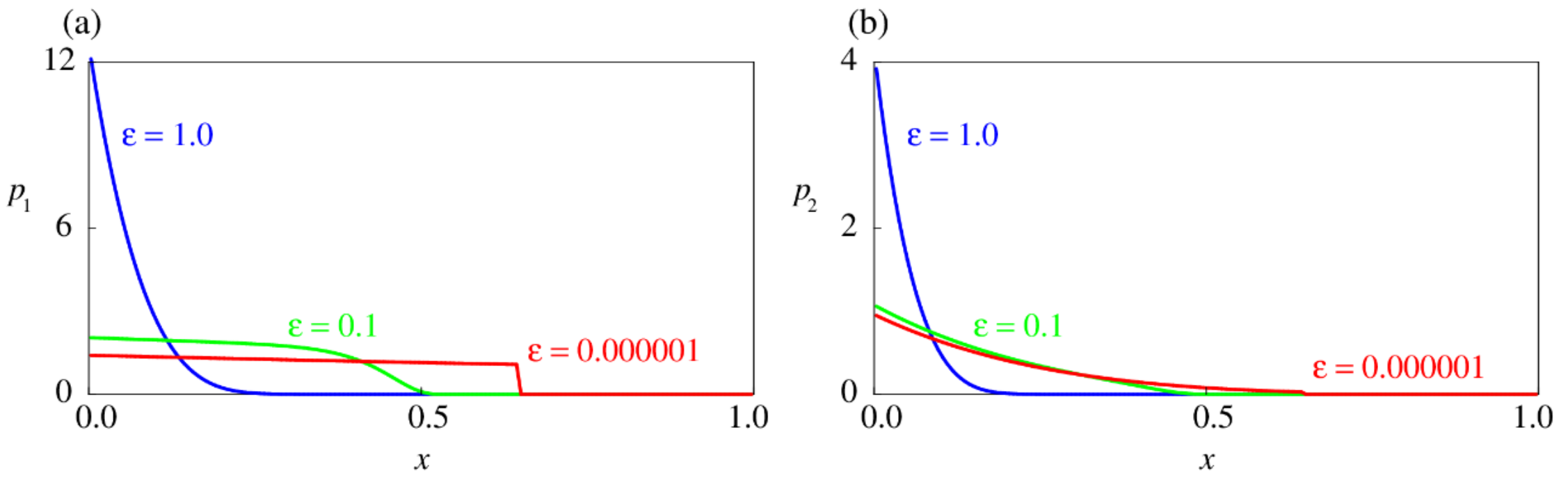

**Fig. 8.** Computed probability densities of MFG with different values of $\varepsilon$ : (a) $i = 1$ and (b) $i = 2$ .

## 5. Conclusion

In this study, a GLD is examined as a special case of PCD in evolutionary game theory and its connection to an MFG of a jump-driven stochastic dynamical system. A generalized logit function in the GLD is investigated such that the proposed theory applies not only to the classical logit dynamic but also covers a wider class of models. The origin of logit functions in the GLD is clarified by their MFG description, showing that the discount rate plays a scaling constant to interpolate both dynamics. The computational analysis focusing on environmental and resource management problems suggests how solutions of the GLD and MFG behave and how the discount rate affect them. After preparing the first version of the manuscript, it has been found that the proposed approach based on an MFG with discount formally works for replicator and projection dynamics [77,78], suggesting its fairly wide applicability.

Player's types are assumed to be fixed in this study, while they can be made time-dependent in problems where group size dynamically changes due to preference changes of players by external factors. The evolution equations in a one-dimensional actions pace are studies, while a space having a higher dimension would arise in applications where players have multiple action variables to be updated. Theoretically, the mathematical and computational approaches employed in this study can be applied to these cases with suitable modifications, but computational cost to discretize the non-local part in the PIDEs can become a bottleneck in case of their numerical implementation. An efficient and stable computational architecture is necessary to resolve the issues. Machine learning schemes can serve as possible candidates; however, key structures of solutions to the PIDEs, such as the conservation and non-negativity of probability and bounds of value functions, require preservation during computation. Other issues to be addressed in the future, which is from an application perspective, involve identifying the utility and cost of the evolutionary game and MFG models from existing social data, and analyzing their regularity because it can affect players' actions both qualitatively and quantitatively. Weaker alternatives to **Assumptions 1 and 2** should also be investigated in the future to better clarify applicability and limitation of the proposed approach.

**Appendix**

***Proof of Proposition 1***

We show the proof for $q \neq 1$. The proof for $q = 1$ is essentially the same because our Utility is bounded by **Assumption 1**. The strategy of proof is basically to find several regularity results on the coefficients of an evolutionary game model (e.g., ), and then apply the generalized Picard–Lindelöf theorem (p. 79 [52]).

First, the three technical lemmas for its proof are discussed. For any $A \in \mathcal{B}$, $\{\mu\} \in (\mathcal{M}_2)^I$ and $i \in \Xi$, set

$$\mathfrak{L}_i \mu(A) = \int_{y \in A} \int_{w \in \Omega} \left\{ \frac{\exp_q\left(\eta^{-1} \Delta U_i\left(y, w, \{\mu\}\right)\right)}{\int_\Omega \exp_q\left(\eta^{-1} \Delta U_i\left(z, w, \{\mu\}\right)\right) \mathrm{d}z} \right\}^q \mathrm{d}w \mathrm{d}y \tag{60}$$

and

$$\bar{\mathfrak{L}}_i \mu(A) = \left(2 - m_i^{-1} \|\mu_i\|\right)_+ \mathfrak{L}_i \mu(t, A). \tag{61}$$

***Lemma 1***

*For any* $A \in \mathcal{B}$, $\{\mu\} \in (\mathcal{M}_2)^I$ *and* $i \in \Xi$, *it follows that*

$$\|\mathfrak{L}_i \mu\|, \theta_i(q, \eta) \leq \left\{ \frac{\exp_q\left(2\eta^{-1} L\right)}{\exp_q\left(-2\eta^{-1} L\right)} \right\}^q \left(= \bar{\theta} \text{ in the sequel}\right). \tag{62}$$

***Proof of Lemma 1***

Note that $\exp_q\left(-2\eta^{-1} L\right) > 0$ by **Assumption 2**. It follows that ("$\mathrm{d}x$" is a one-dimensional Lebesgue measure)

$$\begin{aligned} \|\mathfrak{L}_i \mu\| &= \left\| \int_{w \in \Omega} \left\{ \frac{\exp_q\left(\eta^{-1} \Delta U_i\left(x, w, \{\mu\}\right)\right)}{\int_\Omega \exp_q\left(\eta^{-1} \Delta U_i\left(z, w, \{\mu\}\right)\right) \mathrm{d}z} \right\}^q \mathrm{d}w \mathrm{d}x \right\| \\ &\leq \max_{x, w \in \Omega} \left\{ \frac{\exp_q\left(\eta^{-1} \Delta U_i\left(x, w, \{\mu\}\right)\right)}{\int_\Omega \exp_q\left(\eta^{-1} \Delta U_i\left(z, w, \{\mu\}\right)\right) \mathrm{d}z} \right\}^q \\ &\leq \max_{x, w \in \Omega} \left\{ \frac{\exp_q\left(2\eta^{-1} L\right)}{\int_\Omega \exp_q\left(-2\eta^{-1} L\right) \mathrm{d}z} \right\}^q \\ &= \left\{ \frac{\exp_q\left(2\eta^{-1} L\right)}{\exp_q\left(-2\eta^{-1} L\right)} \right\}^q \end{aligned}. \tag{63}$$

The other inequality follows analogously as

$$
\begin{aligned}
\theta_i(q,\eta) &= \int_{y\in\Omega}\int_{w\in\Omega}\left\{\frac{\exp_q\left(\eta^{-1}\Delta U_i\left(w,y,\{\mu\}\right)\right)}{\int_\Omega \exp_q\left(\eta^{-1}\Delta U_i\left(z,y,\{\mu\}\right)\right)\mathrm{d}z}\right\}^q \mathrm{d}w\mathrm{d}y \\
&\le \max_{y,w\in\Omega}\left\{\frac{\exp_q\left(\eta^{-1}\Delta U_i\left(w,y,\{\mu\}\right)\right)}{\int_\Omega \exp_q\left(\eta^{-1}\Delta U_i\left(z,y,\{\mu\}\right)\right)\mathrm{d}z}\right\}^q \\
&\le \left\{\frac{\exp_q\left(2\eta^{-1}L\right)}{\exp_q\left(-2\eta^{-1}L\right)}\right\}^q
\end{aligned}
\quad . \tag{64}
$$

□

***Lemma 2***

*For any* $\{\mu\},\{\nu\}\in(\mathcal{M}_2)^I$ *and* $i\in\Xi$, *there exists a constant* $l>0$ *such that*

$$
\left\|\mathfrak{L}_i\mu-\mathfrak{L}_i\nu\right\| \le l\sum_{j=1}^{I}\left\|\mu_j-\nu_j\right\|. \tag{65}
$$

***Proof of Lemma 2***

By **Assumption 2**, there is a constant $C_1>0$ independent from $\{\mu\}\in(\mathcal{M}_2)^I$ and $x,z\in\Omega$ such that

$$
\exp_q\left(\eta^{-1}\Delta U_i\left(z,x,\{\mu\}\right)\right)\ge C_1^{-1}. \tag{66}
$$

On fixing some $i\in\Xi$, it follows that

$$
\begin{aligned}
&\left\|\mathfrak{L}_i\mu-\mathfrak{L}_i\nu\right\| \\
&= \left\|\left[\int_{y\in\Omega}\left\{\frac{\exp_q\left(\eta^{-1}\Delta U_i\left(x,y,\{\mu\}\right)\right)}{\int_\Omega \exp_q\left(\eta^{-1}\Delta U_i\left(z,y,\{\mu\}\right)\right)\mathrm{d}z}\right\}^q\mathrm{d}y-\int_{y\in\Omega}\left\{\frac{\exp_q\left(\eta^{-1}\Delta U_i\left(x,y,\{\nu\}\right)\right)}{\int_\Omega \exp_q\left(\eta^{-1}\Delta U_i\left(z,y,\{\nu\}\right)\right)\mathrm{d}z}\right\}^q\mathrm{d}y\right]\mathrm{d}x\right\| \\
&= \left\|\int_{y\in\Omega}\left[\left\{\frac{\exp_q\left(\eta^{-1}\Delta U_i\left(x,y,\{\mu\}\right)\right)}{\int_\Omega \exp_q\left(\eta^{-1}\Delta U_i\left(z,y,\{\mu\}\right)\right)\mathrm{d}z}\right\}^q-\left\{\frac{\exp_q\left(\eta^{-1}\Delta U_i\left(x,y,\{\nu\}\right)\right)}{\int_\Omega \exp_q\left(\eta^{-1}\Delta U_i\left(z,y,\{\nu\}\right)\right)\mathrm{d}z}\right\}^q\right]\mathrm{d}y\mathrm{d}x\right\| \\
&\le \left\|\int_{y\in\Omega}\left|\left\{\frac{\exp_q\left(\eta^{-1}\Delta U_i\left(x,y,\{\mu\}\right)\right)}{\int_\Omega \exp_q\left(\eta^{-1}\Delta U_i\left(z,y,\{\mu\}\right)\right)\mathrm{d}z}\right\}^q-\left\{\frac{\exp_q\left(\eta^{-1}\Delta U_i\left(x,y,\{\nu\}\right)\right)}{\int_\Omega \exp_q\left(\eta^{-1}\Delta U_i\left(z,y,\{\nu\}\right)\right)\mathrm{d}z}\right\}^q\right|\mathrm{d}y\mathrm{d}x\right\| \\
&\le q\max_{\substack{x,y\in\Omega\\ \rho=\mu,\nu}}\left\{\frac{\exp_q\left(\eta^{-1}\Delta U_i\left(x,y,\{\rho\}\right)\right)}{\int_\Omega \exp_q\left(\eta^{-1}\Delta U_i\left(z,y,\{\rho\}\right)\right)\mathrm{d}z}\right\}^{q-1} \\
&\times\left\|\int_{y\in\Omega}\left|\frac{\exp_q\left(\eta^{-1}\Delta U_i\left(x,y,\{\mu\}\right)\right)}{\int_\Omega \exp_q\left(\eta^{-1}\Delta U_i\left(z,y,\{\mu\}\right)\right)\mathrm{d}z}-\frac{\exp_q\left(\eta^{-1}\Delta U_i\left(x,y,\{\nu\}\right)\right)}{\int_\Omega \exp_q\left(\eta^{-1}\Delta U_i\left(z,y,\{\nu\}\right)\right)\mathrm{d}z}\right|\mathrm{d}y\mathrm{d}x\right\| \\
&\le C_2\left\|\int_{y\in\Omega}\left|\begin{aligned}&\exp_q\left(\eta^{-1}\Delta U_i\left(x,y,\{\mu\}\right)\right)\int_\Omega \exp_q\left(\eta^{-1}\Delta U_i\left(z,y,\{\nu\}\right)\right)\mathrm{d}z\\ &-\exp_q\left(\eta^{-1}\Delta U_i\left(x,y,\{\nu\}\right)\right)\int_\Omega \exp_q\left(\eta^{-1}\Delta U_i\left(z,y,\{\mu\}\right)\right)\mathrm{d}z\end{aligned}\right|\mathrm{d}y\mathrm{d}x\right\|
\end{aligned}
\quad . \tag{67}
$$

Herein, $C_2 > 0$ is a constant independent from $x, y \in \Omega$ and $\{\mu\},\{\nu\} \in (\mathcal{M}_2)^I$. **Assumption 1** and (66) lead to

$$q \max_{\substack{x,y\in\Omega \\ \rho=\mu,\nu}} \left\{ \frac{\exp_q\left(\eta^{-1}\Delta U_i\left(x,y,\{\rho\}\right)\right)}{\int_\Omega \exp_q\left(\eta^{-1}\Delta U_i\left(z,y,\{\rho\}\right)\right)\mathrm{dz}} \right\}^{q-1} \begin{aligned} &\leq q \max_{\substack{x,y\in\Omega \\ \rho=\mu,\nu}} \left\{ C_1 \exp_q\left(\eta^{-1}\Delta U_i\left(x,y,\{\rho\}\right)\right)\right\}^{q-1} \\ &\leq q\left\{C_1 \exp_q\left(2\eta^{-1}L\right)\right\}^{q-1} \end{aligned} \quad \text{if } q \geq 1 \tag{68}$$

and

$$\begin{aligned} &q \max_{\substack{x,y\in\Omega \\ \rho=\mu,\nu}} \left\{ \frac{\exp_q\left(\eta^{-1}\Delta U_i\left(x,y,\{\rho\}\right)\right)}{\int_\Omega \exp_q\left(\eta^{-1}\Delta U_i\left(z,y,\{\rho\}\right)\right)\mathrm{dz}} \right\}^{q-1} \\ &= q \max_{\substack{x,y\in\Omega \\ \rho=\mu,\nu}} \left\{ \frac{\int_\Omega \exp_q\left(\eta^{-1}\Delta U_i\left(z,y,\{\rho\}\right)\right)\mathrm{dz}}{\exp_q\left(\eta^{-1}\Delta U_i\left(x,y,\{\rho\}\right)\right)} \right\}^{1-q} \quad \text{if } 0 < q < 1 . \\ &\leq q \max_{\substack{y\in\Omega \\ \rho=\mu,\nu}} \left\{ C_1 \int_\Omega \exp_q\left(\eta^{-1}\Delta U_i\left(z,y,\{\rho\}\right)\right)\mathrm{dz}\right\}^{1-q} \\ &\leq q\left\{C_1 \exp_q\left(2\eta^{-1}L\right)\right\}^{1-q} \end{aligned} \tag{69}$$

Consequently, it follows that

$$q \max_{\substack{x,y\in\Omega \\ \rho=\mu,\nu}} \left\{ \frac{\exp_q\left(\eta^{-1}\Delta U_i\left(x,y,\{\rho\}\right)\right)}{\int_\Omega \exp_q\left(\eta^{-1}\Delta U_i\left(z,y,\{\rho\}\right)\right)\mathrm{dz}} \right\}^{q-1} \leq q\left\{C_1 \exp_q\left(2\eta^{-1}L\right)\right\}^{|q-1|} \tag{70}$$

whose right-hand side is independent from $x, y$, $\{\mu\},\{\nu\}$. By (66), for any $y \in \Omega$, it follows that

$$\int_\Omega \exp_q\left(\eta^{-1}\Delta U_i\left(z,y,\{\mu\}\right)\right)\mathrm{dz} \int_\Omega \exp_q\left(\eta^{-1}\Delta U_i\left(z,y,\{\nu\}\right)\right)\mathrm{dz} \leq C_1^{-2} . \tag{71}$$

We continue as

$$
\begin{aligned}
&\left\|\int_{y\in\Omega}\left|\begin{matrix}\exp_q\left(\eta^{-1}\Delta U_i\left(x,y,\{\mu\}\right)\right)\int_\Omega \exp_q\left(\eta^{-1}\Delta U_i\left(z,y,\{\nu\}\right)\right)\mathrm{d}z\\ -\exp_q\left(\eta^{-1}\Delta U_i\left(x,y,\{\nu\}\right)\right)\int_\Omega \exp_q\left(\eta^{-1}\Delta U_i\left(z,y,\{\mu\}\right)\right)\mathrm{d}z\end{matrix}\right|\mathrm{d}y\mathrm{d}x\right\| \\
&\leq\left\|\int_{y\in\Omega}\left|\begin{matrix}\exp_q\left(\eta^{-1}\Delta U_i\left(x,y,\{\mu\}\right)\right)\int_\Omega \exp_q\left(\eta^{-1}\Delta U_i\left(z,y,\{\nu\}\right)\right)\mathrm{d}z\\ -\exp_q\left(\eta^{-1}\Delta U_i\left(x,y,\{\mu\}\right)\right)\int_\Omega \exp_q\left(\eta^{-1}\Delta U_i\left(z,y,\{\mu\}\right)\right)\mathrm{d}z\end{matrix}\right|\mathrm{d}y\mathrm{d}x\right\| \\
&+\left\|\int_{y\in\Omega}\left|\begin{matrix}\exp_q\left(\eta^{-1}\Delta U_i\left(x,y,\{\mu\}\right)\right)\int_\Omega \exp_q\left(\eta^{-1}\Delta U_i\left(z,y,\{\mu\}\right)\right)\mathrm{d}z\\ -\exp_q\left(\eta^{-1}\Delta U_i\left(x,y,\{\nu\}\right)\right)\int_\Omega \exp_q\left(\eta^{-1}\Delta U_i\left(z,y,\{\mu\}\right)\right)\mathrm{d}z\end{matrix}\right|\mathrm{d}y\mathrm{d}x\right\| \\
&=\left\|\int_{y\in\Omega}\left|\begin{matrix}\exp_q\left(\eta^{-1}\Delta U_i\left(x,y,\{\mu\}\right)\right)\times\\ \int_\Omega\left\{\exp_q\left(\eta^{-1}\Delta U_i\left(z,y,\{\nu\}\right)\right)-\exp_q\left(\eta^{-1}\Delta U_i\left(z,y,\{\mu\}\right)\right)\right\}\mathrm{d}z\end{matrix}\right|\mathrm{d}y\mathrm{d}x\right\| \\
&+\left\|\int_{y\in\Omega}\left|\begin{matrix}\int_\Omega \exp_q\left(\eta^{-1}\Delta U_i\left(z,y,\{\mu\}\right)\right)\mathrm{d}z\times\\ \left\{\exp_q\left(\eta^{-1}\Delta U_i\left(x,y,\{\mu\}\right)\right)-\exp_q\left(\eta^{-1}\Delta U_i\left(x,y,\{\nu\}\right)\right)\right\}\end{matrix}\right|\mathrm{d}y\mathrm{d}x\right\| \\
&\leq \exp_q\left(2\eta^{-1}L\right)\left\|\int_{y\in\Omega}\left|\int_{z\in\Omega}\left\{\begin{matrix}\exp_q\left(\eta^{-1}\Delta U_i\left(z,y,\{\nu\}\right)\right)\\ -\exp_q\left(\eta^{-1}\Delta U_i\left(z,y,\{\mu\}\right)\right)\end{matrix}\right\}\mathrm{d}z\right|\mathrm{d}y\mathrm{d}x\right\| \\
&+\exp_q\left(2\eta^{-1}L\right)\left\|\int_{y\in\Omega}\left|\exp_q\left(\eta^{-1}\Delta U_i\left(x,y,\{\mu\}\right)\right)-\exp_q\left(\eta^{-1}\Delta U_i\left(x,y,\{\nu\}\right)\right)\right|\mathrm{d}y\mathrm{d}x\right\| \\
&\leq 2\exp_q\left(2\eta^{-1}L\right)\left\|\int_{y\in\Omega}\left|\exp_q\left(\eta^{-1}\Delta U_i\left(x,y,\{\mu\}\right)\right)-\exp_q\left(\eta^{-1}\Delta U_i\left(x,y,\{\nu\}\right)\right)\right|\mathrm{d}y\mathrm{d}x\right\|
\end{aligned}
\quad . \qquad (72)
$$

Here, $0\leq \exp_q\left(\eta^{-1}\Delta U_i\left(x,z,\{\mu\}\right)\right)\leq \exp_q\left(2\eta^{-1}L\right)$ for any $x,z\in\Omega$ and $\{\mu\}\in\left(\mathcal{M}_2\right)^I$ by **Assumption 1**. For any $x,y\in\Omega$, it follows that

$$
\begin{aligned}
&\left|\exp_q\left(\eta^{-1}\Delta U_i\left(x,y,\{\mu\}\right)\right)-\exp_q\left(\eta^{-1}\Delta U_i\left(x,y,\{\nu\}\right)\right)\right| \\
&\leq\left\{\exp_q\left(2\eta^{-1}L\right)\right\}^q\left\{\left|U_i\left(x,\{\mu\}\right)-U_i\left(x,\{\nu\}\right)\right|+\left|U_i\left(y,\{\mu\}\right)-U_i\left(y,\{\nu\}\right)\right|\right\} \\
&\leq 2\left\{\exp_q\left(2\eta^{-1}L\right)\right\}^q\max\left\{\left|U_i\left(x,\{\mu\}\right)-U_i\left(x,\{\nu\}\right)\right|,\left|U_i\left(y,\{\mu\}\right)-U_i\left(y,\{\nu\}\right)\right|\right\} \\
&\leq 2L\left\{\exp_q\left(2\eta^{-1}L\right)\right\}^q\sum_{j=1}^{I}\left\|\mu_j-\nu_j\right\|
\end{aligned}
\quad . \qquad (73)
$$

Consequently, by (67), (72), (73), $\int_{y\in\Omega}\mathrm{d}y=1$, and $\left\|\mathrm{d}x\right\|=1$, we get (65) with a sufficiently large $l>0$.

□

***Lemma 3***

*For any* $\{\mu\},\{\nu\}\in\left(\mathcal{M}\right)^I$ *and* $i\in\Xi$, *there exists a constant* $\bar{l}>0$ *such that*

$$\left\|\bar{\mathfrak{L}}_i\mu-\bar{\mathfrak{L}}_i\nu\right\|\leq\bar{l}\sum_{j=1}^{I}\left\|\mu_j-\nu_j\right\| . \qquad (74)$$

***Proof of Lemma 3***

By fixing some $i \in \Xi$, the four cases: $2m_i \le \|\mu_i\|$ and $2m_i \le \|\nu_i\|$; $2m_i \ge \|\mu_i\|$ and $2m_i \ge \|\nu_i\|$; $2m_i \ge \|\mu_i\|$ and $2m_i \le \|\nu_i\|$; $2m_i \le \|\mu_i\|$ and $2m_i \ge \|\nu_i\|$, are considered.

For the first case, nothing is left to prove as $\left(2 - m_i^{-1}\|\mu_i(t,\cdot)\|\right)_+ = \left(2 - m_i^{-1}\|\nu_i(t,\cdot)\|\right)_+ = 0$. For the second case, we have

$$\begin{aligned} \left\|\bar{\mathfrak{L}}_i\mu - \bar{\mathfrak{L}}_i\nu\right\| &= \left\|\left(2 - m_i^{-1}\|\mu_i\|\right)\mathfrak{L}_i\mu - \left(2 - m_i^{-1}\|\nu_i\|\right)\mathfrak{L}_i\nu\right\| \\ &\le 2\|\mathfrak{L}_i\mu - \mathfrak{L}_i\nu\| + m_i^{-1}\left\|\|\mu_i\|\mathfrak{L}_i\mu - \|\nu_i\|\mathfrak{L}_i\nu\right\| \\ &\le 2\|\mathfrak{L}_i\mu - \mathfrak{L}_i\nu\| + m_i^{-1}\left\|\|\mu_i\|\mathfrak{L}_i\mu - \|\mu_i\|\mathfrak{L}_i\nu\right\| + m_i^{-1}\left|\|\mu_i\|\mathfrak{L}_i\nu - \|\nu_i\|\mathfrak{L}_i\nu\right| . \\ &\le 2\|\mathfrak{L}_i\mu - \mathfrak{L}_i\nu\| + m_i^{-1}\|\mu_i\|\|\mathfrak{L}_i\mu - \mathfrak{L}_i\nu\| + m_i^{-1}\|\mathfrak{L}_i\nu\|\left|\|\mu_i\| - \|\nu_i\|\right| \\ &\le \left(2 + m_i^{-1}\|\mu_i\|\right)\|\mathfrak{L}_i\mu - \mathfrak{L}_i\nu\| + m_i^{-1}\|\mathfrak{L}_i\nu\|\left|\|\mu_i\| - \|\nu_i\|\right| \end{aligned} \tag{75}$$

Then, it follows that

$$2 + m_i^{-1}\|\mu_i\| \le 2 + m_i^{-1} \times 2m_i = 4 \quad \text{and} \quad \left|\|\mu_i\| - \|\nu_i\|\right| \le \|\mu_i - \nu_i\|, \tag{76}$$

and hence the inequality

$$\begin{aligned} \left\|\bar{\mathfrak{L}}_i\mu - \bar{\mathfrak{L}}_i\nu\right\| &\le 4\|\mathfrak{L}_i\mu - \mathfrak{L}_i\nu\| + m_i^{-1}\|\mathfrak{L}_i\nu\|\|\mu_i - \nu_i\| \\ &\le 4l\|\mu_i - \nu_i\| + m_i^{-1}\bar{\theta}\|\mu_i - \nu_i\| \qquad , \\ &\le \left(4l + m_i^{-1}\bar{\theta}\right)\sum_{j=1}^{I}\|\mu_j - \nu_j\| \end{aligned} \tag{77}$$

where $\|\mathfrak{L}_i\nu\| \le \bar{\theta}$ by **Lemma 1** due to $\nu_i \in \mathcal{P}_{2m_i} \in \mathcal{M}_2$, and **Lemma 2** was used as well. From (77), it follows that

$$\left\|\bar{\mathfrak{L}}_i\mu - \bar{\mathfrak{L}}_i\nu\right\| \le \left(4l + \bar{\theta}m_i^{-1}\right)\sum_{j=1}^{I}\|\mu_j - \nu_j\| \le \left(4l + \bar{\theta}\max_{j\in\Xi} m_j^{-1}\right)\sum_{j=1}^{I}\|\mu_j - \nu_j\| . \tag{78}$$

For the third case, $\bar{\mathfrak{L}}_i\nu = 0$, and hence

$$\begin{aligned} \left\|\bar{\mathfrak{L}}_i\mu - \bar{\mathfrak{L}}_i\nu\right\| &= \left\|\left(2 - m_i^{-1}\|\mu_i\|\right)\mathfrak{L}_i\mu\right\| \\ &= \|\mathfrak{L}_i\mu\|\left(2 - m_i^{-1}\|\mu_i\|\right) \\ &\le \bar{\theta}\left(2 - m_i^{-1}\|\mu_i\|\right) \\ &\le \bar{\theta}\left(m_i^{-1}\|\nu_i\| - m_i^{-1}\|\mu_i\|\right) . \\ &= \bar{\theta}m_i^{-1}\left(\|\nu_i\| - \|\mu_i\|\right) \\ &\le \bar{\theta}\max_{j\in\Xi} m_j^{-1}\sum_{j=1}^{I}\|\mu_j - \nu_j\| \end{aligned} \tag{79}$$

The fourth case is symmetrical to the third case, and the same $\bar{l}$ is obtained. Consequently, the statement of the proposition with a sufficiently large $\bar{l} > 0$ is proven.

□

By **Lemma 3**, for any $\{\mu\},\{\nu\} \in (\mathcal{M})^I$, $i \in \Xi$, it follows that

$$
\begin{aligned}
\left\| m_i \bar{\mathfrak{L}}_i \mu - \theta_i(q,\eta)\mu_i - \left( m_i \bar{\mathfrak{L}}_i \nu - \theta_i(q,\eta)\nu_i \right) \right\| &\le m_i \left\| \bar{\mathfrak{L}}_i \mu - \bar{\mathfrak{L}}_i \nu \right\| + \theta_i(q,\eta) \left\| \mu_i - \nu_i \right\| \\
&\le \left( \bar{l} + \bar{\theta} \right) \sum_{j=1}^{I} \left\| \mu_j - \nu_j \right\|
\end{aligned}
. \tag{80}
$$

Moreover, for any $\{\mu\} \in (\mathcal{M})^I$ and $i \in \Xi$, the following inequality holds true:

$$
\left\| m_i \bar{\mathfrak{L}}_i \mu - \theta_i(q,\eta)\mu_i \right\| \le m_i \left\| \bar{\mathfrak{L}}_i \mu \right\| + \theta_i(q,\eta) \left\| \mu_i \right\| \le \bar{\theta}\left(1 + \left\| \mu_i \right\|\right) < +\infty . \tag{81}
$$

Considering (80) and (81), the generalized Picard–Lindelöf theorem (Theorem A.3 [41] and p. 79 [52]) applies to the modified GLD (15), showing that it admits a unique solution $\{\mu(t)\} = \{\mu_i(t)\}_{i=1,2,3,\ldots,I} \in (\mathcal{M})^I$ ($t \ge 0$).

Consider the unique solution $\{\mu(t)\} \in (\mathcal{M})^I$ ($t \ge 0$) to the modified GLD (15). For any $A \in \mathcal{B}$, $i \in \Xi$, and $t > 0$, due to $\mathfrak{L}_i \mu(t,A) \ge 0$, it follows that

$$
\begin{aligned}
\frac{\mathrm{d}}{\mathrm{d}t} \mu_i(t,A) &= \left( 2 - m_i^{-1} \left\| \mu_i(t,\cdot) \right\| \right)_+ \left( m_i \mathfrak{L}_i(t,A) - \theta_i(q,\eta)\mu_i(t,A) \right) \\
&\ge -\left( 2 - m_i^{-1} \left\| \mu_i(t,\cdot) \right\| \right)_+ \theta_i(q,\eta)\mu_i(t,A) \\
&\ge -\left( 2 - m_i^{-1} \left\| \mu_i(t,\cdot) \right\| \right)_+ \bar{\theta}\mu_i(t,A)
\end{aligned}
. \tag{82}
$$

When combined with the initial condition $\mu_i(0) \in \mathcal{P}_{m_i}$, $\mu_i(t,A) \ge 0$ is demonstrated, and hence $\mu_i(t,\cdot)$ is a non-negative measure. Additionally, (recall the PCD representation (16))

$$
\frac{\mathrm{d}}{\mathrm{d}t} \mu_i(t,\Omega) = \int_\Omega \frac{\mathrm{d}}{\mathrm{d}t} \mu_i(t,\mathrm{d}x) = \theta_i(q,\eta)\left\{ m_i - \mu_i(t,\Omega) \right\} \quad \text{for any} \quad t > 0 . \tag{83}
$$

Note that

$$
\begin{aligned}
& \int_{x\in\Omega} \int_{y\in\Omega} \rho_i(x,y,\{\mu\}) \mu_i(t,\mathrm{d}y)\mathrm{d}x \\
&= \int_{x\in\Omega} \int_{y\in\Omega} \int_{w\in\Omega} \left\{ \frac{\exp_q\left( \eta^{-1}\Delta U_i(x,w,\{\mu\}) \right)}{\int_\Omega \exp_q\left( \eta^{-1}\Delta U_i(z,w,\{\mu\}) \right)\mathrm{d}z} \right\}^q \mathrm{d}w \mu_i(\mathrm{d}y)\mathrm{d}x \\
&= \int_{x\in\Omega} \int_{w\in\Omega} \left\{ \frac{\exp_q\left( \eta^{-1}\Delta U_i(x,w,\{\mu\}) \right)}{\int_\Omega \exp_q\left( \eta^{-1}\Delta U_i(z,w,\{\mu\}) \right)\mathrm{d}z} \right\}^q \mathrm{d}w\mathrm{d}x \left( \int_{y\in\Omega} \mu_i(\mathrm{d}y) \right) . \\
&= m_i \int_{w\in\Omega} \int_{x\in\Omega} \left\{ \frac{\exp_q\left( \eta^{-1}\Delta U_i(x,w,\{\mu\}) \right)}{\int_\Omega \exp_q\left( \eta^{-1}\Delta U_i(z,w,\{\mu\}) \right)\mathrm{d}z} \right\}^q \mathrm{d}x\mathrm{d}w \\
&= m_i \theta_i(q,\eta)
\end{aligned}
\tag{84}
$$

The ordinary differential equation (83) combined with $\mu_i(0,\Omega) = m_i$ yields $\mu_i(t,\Omega) = m_i$ for any $t > 0$, showing that $\mu_i(t) \in \mathcal{P}_{m_i}$ as desired. Therefore, the unique solution to the modified GLD (15) is a solution to the GLD (6). Finally, the unique solvability of the GLD (6) in the sense of $\mu_i(t) \in \mathcal{P}_{m_i}$ ($t \ge 0$) follows from the Lipschitz continuity estimate in **Lemma 1** combined with a classical Gronwall lemma, which completes the proof of **Proposition 1**.

***Proof of Proposition 2***

The statement in the first sentence of the proposition follows from the method of doubling variables for nonlocal PIDEs (e.g., Theorem 3.4 [79]).

To show the statement in the second sentence, it suffices to show that constant functions $L$ and $L$ serve are viscosity super-solution ad sub-solution, respectively, such that

$$H_{\hat{i}}\left(\hat{t},\hat{x},L,L,\frac{\partial\bar{\Psi}^{(\hat{i})}\left(\hat{t},\hat{x}\right)}{\partial t}\right)\geq 0 \quad\text{and}\quad H_{\hat{i}}\left(\hat{t},\hat{x},-L,-L,\frac{\partial\underline{\Psi}^{(\hat{i})}\left(\hat{t},\hat{x}\right)}{\partial t}\right)\leq 0 \tag{85}$$

for proper test functions $\bar{\Psi}^{(\hat{i})}$ and $\underline{\Psi}^{(\hat{i})}$, where the other notations are the same with those in **Definition 1**. The left inequality in (85) is equivalent to

$$-\frac{\partial\bar{\Psi}^{(\hat{i})}\left(\hat{t},\hat{x}\right)}{\partial t}+\delta L-\delta U_{\hat{i}}\left(\hat{x},\left\{\mu\left(\hat{t},\cdot\right)\right\}\right)\geq 0, \tag{86}$$

but this is satisfied because of $\left|U_i\right|\leq L$ and test functions satisfy $\frac{\partial\bar{\Psi}^{(\hat{i})}\left(\hat{t},\hat{x}\right)}{\partial t}\leq 0$. Herein, any smooth test function $\bar{\Psi}^{(\hat{i})}$ touches the constant function $L$ from below at any points on $[0,T)\times\Omega$, and hence its partial differential $\frac{\partial\bar{\Psi}^{(\hat{i})}\left(\hat{t},\hat{x}\right)}{\partial t}$ vanishes if $\hat{t}\in(0,T)$, and is non-negative if $\hat{t}=0$. The other term in the equality portion of (86) is proven in a similar manner, thereby concluding the proof.

***Proof of Proposition 3***

An induction argument starting from $k=N_t$ is considered. The inequality (52) is trivially satisfied at $k=N_t$. Assuming that the inequality is satisfied at some $k\in\{1,2,3,...,N_t\}$, by (43) it follows that

$$\Phi_{k-1,l}^{(i)}=\Phi_{k,l}^{(i)}+\Delta t\times\left\{\eta\ln_q\left(\sum_{m=1}^{I}\exp_q\left(\eta^{-1}\left(\Phi_{k,m}^{(i)}-\Phi_{k,l}^{(i)}\right)\right)\Delta x\right)-\delta\Phi_{k,l}^{(i)}+\delta U_{i,k-1,l}\right\}, \tag{87}$$

where the right-hand side terms are well-defined as a real value due to **Assumption 2** and the inequality (52) for $\Phi_{k,\cdot}^{(i)}$. The classical Jensen's lemma gives

$$\eta\ln_q\left(\sum_{m=1}^{I}\exp_q\left(\eta^{-1}\left(\Phi_{k,m}^{(i)}-\Phi_{k,l}^{(i)}\right)\right)\Delta x\right)\geq\eta\sum_{m=1}^{I}\eta^{-1}\left(\Phi_{k,m}^{(i)}-\Phi_{k,l}^{(i)}\right)\Delta x=\sum_{m=1}^{I}\left(\Phi_{k,m}^{(i)}-\Phi_{k,l}^{(i)}\right)\Delta x, \tag{88}$$

and hence

$$\begin{aligned}
\Phi_{k-1,l}^{(i)} &\geq \Phi_{k,l}^{(i)} + \Delta t \times \left\{ \sum_{m=1}^{I} \left( \Phi_{k,m}^{(i)} - \Phi_{k,l}^{(i)} \right) \Delta x - \delta \Phi_{k,l}^{(i)} + \delta U_{i,k-1,l} \right\} \\
&\geq \Phi_{k,l}^{(i)} + \Delta t \times \left\{ -\sum_{m=1}^{I} \Phi_{k,l}^{(i)} \Delta x - \delta \Phi_{k,l}^{(i)} \right\} \\
&\geq \Phi_{k,l}^{(i)} \left\{ 1 - \Delta t - \delta \Delta t \right\} \\
&\geq 0
\end{aligned} \tag{89}$$

due to the non-negativity of $\Phi_{k,\cdot}^{(i)}$, $U_{i,\cdot,\cdot}$, and the assumption that $\Delta t \leq \dfrac{1}{1+\delta}$, where we used $\sum_{m=1}^{I} \Delta x = 1$.

The other side of the inequality (52) is proven by assuming the induction argument. If the maximum value among $\Phi_{k,l}^{(i)}$ ($l = 1,2,3,...,N_x$) is attained at some $l = l'$, then $\Phi_{k,l}^{(i)} \leq L$, and

$$\begin{aligned}
\Phi_{k-1,l}^{(i)} &\leq \Phi_{k,l}^{(i)} + \Delta t \times \left\{ \eta \ln_q \left( \sum_{m=1}^{I} \exp_q \left( \eta^{-1} \left( \Phi_{k,m}^{(i)} - \Phi_{k,l}^{(i)} \right) \right) \Delta x \right) - \delta \Phi_{k,l}^{(i)} + \delta U_{i,k-1,l} \right\} \\
&\leq \Phi_{k,l}^{(i)} \left( 1 - \delta \Delta t \right) + \Delta t \eta \ln_q \left( \sum_{m=1}^{I} \exp_q \left( \eta^{-1} \left( L - \Phi_{k,l}^{(i)} \right) \right) \Delta x \right) + \delta \Delta t L \\
&= \Phi_{k,l}^{(i)} \left( 1 - \delta \Delta t \right) + \Delta t \eta \ln_q \left( \exp_q \left( \eta^{-1} \left( L - \Phi_{k,l}^{(i)} \right) \right) \right) + \delta \Delta t L \\
&= \Phi_{k,l}^{(i)} \left( 1 - \delta \Delta t \right) + \Delta t \eta \eta^{-1} \left( L - \Phi_{k,l}^{(i)} \right) + \delta \Delta t L \\
&= \Phi_{k,l}^{(i)} \left( 1 - \Delta t - \delta \Delta t \right) + \left( 1 + \delta \right) \Delta t L \\
&\leq L \left( 1 - \Delta t - \delta \Delta t \right) + \left( 1 + \delta \right) \Delta t L \\
&= L
\end{aligned} , \tag{90}$$

with which the proof is completed by the induction.

***Proof of Proposition 4***

An induction argument starting from $k = 0$. The inequalities (53) and (54) are trivially satisfied at $k = 0$. Assuming that the inequality is satisfied at some $k \in \{0,1,2,...,N_t - 1\}$, using (44) for each $i \in \Xi$ and $m,l = 1,2,3,...,N_x$, it follows that

$$\varphi_{i,k,m,l}^{*} = \frac{\exp_q \left( \eta^{-1} \left( \Phi_{k,m}^{(i)} - \Phi_{k,l}^{(i)} \right) \right)}{\sum_{o=1}^{I} \exp_q \left( \eta^{-1} \left( \Phi_{k,o}^{(i)} - \Phi_{k,l}^{(i)} \right) \right) \Delta x} \geq 0 \tag{91}$$

and

$$\varphi_{i,k,m,l}^{*} = \frac{\exp_q \left( \eta^{-1} \left( \Phi_{k,m}^{(i)} - \Phi_{k,l}^{(i)} \right) \right)}{\sum_{o=1}^{I} \exp_q \left( \eta^{-1} \left( \Phi_{k,o}^{(i)} - \Phi_{k,l}^{(i)} \right) \right) \Delta x} \leq \frac{\exp_q \left( 2\eta^{-1} L \right)}{\sum_{o=1}^{I} \exp_q \left( -2\eta^{-1} L \right) \Delta x} = \frac{\exp_q \left( 2\eta^{-1} L \right)}{\exp_q \left( -2\eta^{-1} L \right)} . \tag{92}$$

Then, the following inequality is obtained:

$$
\begin{aligned}
\mu_{i,k+1,l} &\geq \mu_{i,k,l} + \Delta t \left\{ \sum_{m=1}^{N_x} \left( \varphi_{i,k,l,m}^{*} \right)^{q} \Delta x - \left( \sum_{m=1}^{N_x} \left( \frac{\exp_q\left(2\eta^{-1}L\right)}{\exp_q\left(-2\eta^{-1}L\right)} \right)^{q} \Delta x \right) \mu_{i,k,l} \right\} \\
&\geq \mu_{i,k,l} + \Delta t \left\{ -\left( \frac{\exp_q\left(2\eta^{-1}L\right)}{\exp_q\left(-2\eta^{-1}L\right)} \right)^{q} \mu_{i,k,l} \right\} \\
&\geq \mu_{i,k,l} \left\{ 1 - \left( \frac{\exp_q\left(2\eta^{-1}L\right)}{\exp_q\left(-2\eta^{-1}L\right)} \right)^{q} \Delta t \right\} \\
&\geq 0
\end{aligned}
\quad . \tag{93}
$$

Moreover, the following equality holds true with which the proof can be completed by induction:

$$
\begin{aligned}
\sum_{l=1}^{N_x} \mu_{i,k+1,l} &= \sum_{l=1}^{N_x} \mu_{i,k,l} + \Delta t \sum_{l=1}^{N_x} \left\{ \sum_{m=1}^{N_x} \left( \varphi_{i,k,l,m}^{*} \right)^{q} \mu_{i,k,m} \Delta x - \left( \sum_{m=1}^{N_x} \left( \varphi_{i,k,m,l}^{*} \right)^{q} \Delta x \right) \mu_{i,k,l} \right\} \\
&= m_i + \Delta t \left\{ \sum_{l,m=1}^{N_x} \left( \varphi_{i,k,l,m}^{*} \right)^{q} \mu_{i,k,m} \Delta x - \sum_{l,m=1}^{N_x} \left( \varphi_{i,k,m,l}^{*} \right)^{q} \mu_{i,k,l} \Delta x \right\} \\
&= m_i + \Delta t \left\{ \sum_{l,m=1}^{N_x} \left( \varphi_{i,k,l,m}^{*} \right)^{q} \mu_{i,k,m} \Delta x - \sum_{l,m=1}^{N_x} \left( \varphi_{i,k,l,m}^{*} \right)^{q} \mu_{i,k,m} \Delta x \right\} \\
&= m_i
\end{aligned}
\quad . \tag{94}
$$

**Competing interests:** The author declares no competing interest.

**Funding:** Japan Society for the Promotion of Science 22K14441.

**Availability of data and material:** Data are available upon reasonable request to the corresponding author.

**Declaration:** No AI technologies were used for writing the manuscript.

**Author contributions:**

*Hidekazu Yoshioka*: Conceptualization, Methodology, Software, Formal analysis, Data curation, Visualization, Writing–original draft preparation, Writing–review and editing, Supervision, Project administration, Funding acquisition.